\documentclass[thmsa,11pt]{article}
\usepackage{amsmath}
\usepackage{amssymb}

\begin{document}

\title{Fixed Points of Self-embeddings of Models of Arithmetic}
\author{Saeideh Bahrami \& Ali Enayat}
\maketitle

\begin{abstract}
\noindent We investigate the structure of \textit{fixed point sets }of
self-embeddings of models of arithmetic. Our principal results are Theorems
A, B, and C below. \medskip

\noindent In what follows $\mathcal{M}$ is a countable nonstandard model of
the fragment $\mathrm{I}\Sigma _{1}$ of \textrm{PA} (Peano Arithmetic); $%
\mathbb{N}$ is the initial segment of $\mathcal{M}$ consisting of standard
numbers of $\mathcal{M}$; $\mathrm{I}_{\mathrm{fix}}(j)$ is the longest
initial segment of fixed points of $j$; $\mathrm{Fix}(j)$ is the fixed point
set of $j$; $K^{1}(\mathcal{M})$ consists of $\Sigma _{1}$-definable
elements of $\mathcal{M}$; and a self-embedding $j$ of $\mathcal{M}$ is said
to be a proper initial self-embedding if $j(\mathcal{M})$\ is a proper
initial segment of $\mathcal{M}$.\medskip

\noindent \textbf{Theorem A}.~\textit{The following are equivalent for a
proper initial segment }$I$ \textit{of} $\mathcal{M}$:

\noindent \textbf{(1)} $I=\mathrm{I}_{\mathrm{fix}}(j)$ \textit{for some
self-embedding }$j$ \textit{of} $\mathcal{M}$.

\noindent \textbf{(2)} $I$\textit{\ is closed under exponentiation.}

\noindent \textbf{(3)} $I=$ $\mathrm{I}_{\mathrm{fix}}(j)$ \textit{for some
proper initial self-embedding }$j$ \textit{of} $\mathcal{M}$.\medskip

\noindent \textbf{Theorem B}.~\textit{The following are equivalent for a
proper initial segment }$I$ \textit{of} $\mathcal{M}$:

\noindent \textbf{(1)} $I=$ $\mathrm{Fix}(j)$ \textit{for some
self-embedding }$j$ \textit{of} $\mathcal{M}$.

\noindent \textbf{(2)} $I$\textit{\ is a strong cut of }$\mathcal{M}$\textit{%
\ and }$I\prec _{\Sigma _{1}}\mathcal{M}$.

\noindent \textbf{(3)} $I=\mathrm{Fix}(j)$ \textit{for some proper initial
self-embedding }$j$ \textit{of} $\mathcal{M}$.\medskip

\noindent \textbf{Theorem C}.~\textit{The following are equivalent}:

\noindent \textbf{(1)} $\mathrm{Fix}(j)=K^{1}(\mathcal{M})$ \textit{for some
self-embedding }$j$ \textit{of} $\mathcal{M}$.

\noindent \textbf{(2)} $\mathbb{N}$ \textit{is a strong cut of }$\mathcal{M}$%
.

\noindent \textbf{(3)} $\mathrm{Fix}(j)=K^{1}(\mathcal{M})$ \textit{for some
proper initial self-embedding }$j$ \textit{of} $\mathcal{M}$.

\begin{equation*}
\ast \ast \ast \ast \ast \ast \ast \ast \ast \ast
\end{equation*}

\smallskip

\end{abstract}

\pagebreak

\begin{center}
\textbf{1.}~\textbf{INTRODUCTION\medskip }
\end{center}

In the early 1970s Harvey Friedman \cite[Thm.~4.4]{Harvey} proved a
remarkable theorem: Every countable nonstandard model $\mathcal{M}$ of $%
\mathrm{PA}$ carries a proper initial self-embedding $j$; i.e., $j$
isomorphically maps $\mathcal{M}$ onto a proper initial segment of $\mathcal{%
M}$. Friedman's theorem has been generalized and refined in several ways
over the past several decades (most recently in \cite{Yokoyama} and \cite%
{Ali-Lawrence}). In the mid-1980s Ressayre \cite{Ressayre}, and
independently Dimitracopoulos \& Paris \cite{Costas-Jeff}, generalized
Friedman's theorem by weakening \textrm{PA} to the fragment $\mathrm{I}%
\Sigma _{1}$ of \textrm{PA}. In this paper we refine their work by
investigating \textit{fixed point sets} of self-embeddings of countable
nonstandard models of $\mathrm{I}\Sigma _{1}$.\medskip

Our work here was inspired by certain striking results concerning the
structure of fixed point sets of \textit{automorphisms of countable
recursively saturated models of }$\mathrm{PA}$ summarized in Theorem 1.1
below. In what follows $\mathbb{N}$ is the initial segment of $\mathcal{M}$
consisting of the standard numbers of $\mathcal{M}$;\textit{\ }$K(\mathcal{M}%
)$ is the set of definable elements of $\mathcal{M}$; $\mathrm{I}_{\mathrm{%
fix}}(j)$ is the longest initial segment of fixed points of $j$; and $%
\mathrm{Fix}(j)$ is the fixed point set of $j$, in other words:

\begin{center}
$\mathrm{I}_{\mathrm{fix}}(j):=\{m\in M:\forall x\leq m\ j(x)=x\},$ and $%
\mathrm{Fix}(j):=\{m\in M:j(m)=m\}.$
\end{center}

\noindent \textbf{1.1.}~\textbf{Theorem.}~\textit{Suppose }$\mathcal{M}$%
\textit{\ is a countable recursively saturated model of }$\mathrm{PA}$%
\textit{, and }$I$\textit{\ is a proper initial segment of }$\mathcal{M}$.%
\textit{\ \smallskip }

\noindent \textbf{(a)}\textit{\ }(Smory\'{n}ski \cite{Smorynski-Ifix})%
\textit{\ }$I=$\textit{\ }$\mathrm{I}_{\mathrm{fix}}(j)$ \textit{for} \emph{%
some automorphism }$j$ \emph{of} $\mathcal{M}$ \emph{iff} $I$\textit{\ }%
\emph{is closed under exponentiation}$.$\footnote{%
Smory\'{n}ski established the right-to-left direction of this result and
left the status of the other, much easier direction as an open problem. It
is unclear who first established the easier direction, but by now it is
considered part of the folklore of the subject. A different proof of (a
stronger version of) Smory\'{n}ski's theorem was established in \cite{Me FM}.%
}\textit{\ \smallskip }

\noindent \textbf{(b)}\textit{\ }(Kaye-Kossak-Kotlarski \cite%
{Richard-Roman-Henryk})\textit{\ }$I=\mathrm{Fix}(j)$ \emph{for some
automorphism }$j$ \emph{of} $\mathcal{M}$ \emph{iff}\textit{\ }$(I$ \emph{is
a strong cut of }$\mathcal{M}$ \emph{and} $I\prec \mathcal{M)}$.\textit{%
\smallskip }

\noindent \textbf{(c)}\textit{\ }(Kaye-Kossak-Kotlarski \cite%
{Richard-Roman-Henryk})\textit{\ }$\mathrm{Fix}(j)=$ $K(\mathcal{M})$ \emph{%
for some automorphism }$j$ \emph{of} $\mathcal{M}$ \emph{iff} $\mathbb{N}$
\textit{is a strong cut of }$\mathcal{M}$.\footnote{%
This result was generalized in \cite{Me unified} by showing that if $\mathbb{%
N}$ is strong in $\mathcal{M}$, then the isomorphism types of fixed point
sets of automorphisms of $\mathcal{M}$ are precisely the isomorphism types
of elementary submodels of $\mathcal{M}$, thus confirming a conjecture of
Schmerl.}\textit{\medskip }

In this paper we formulate and establish appropriate analogues of each part
of Theorem 1.1 for \textit{self-embeddings} \textit{of countable nonstandard
models of} $\mathrm{I}\Sigma _{1}$, as encapsulated in Theorem 1.2 below. In
part (c), $K^{1}(\mathcal{M})$ consists of $\Sigma _{1}$-definable elements
of $\mathcal{M}$.\medskip

\noindent \textbf{1.2.}~\textbf{Theorem.}~\textit{Suppose }$\mathcal{M}$%
\textit{\ is a countable nonstandard model of }$\mathrm{I}\Sigma _{1},$
\textit{and }$I$\textit{\ is a proper initial segment of }$\mathcal{M}$%
\textit{. \smallskip }

\noindent \textbf{(a)}\textit{\ }$I=\mathrm{I}_{\mathrm{fix}}(j)$ \emph{for
some self-embedding }$j$ \emph{of} $\mathcal{M}$ \emph{iff} $I$\textrm{\ }%
\emph{is closed under exponentiation }\textit{iff }$I=\mathrm{I}_{\mathrm{fix%
}}(j)$ \emph{for some proper initial self-embedding }$j$ \emph{of} $\mathcal{%
M}$.\textit{\smallskip }

\noindent \textbf{(b)}\textit{\ }$I=\mathrm{Fix}(j)$ \emph{for some} \emph{%
self-embedding }$j$ \emph{of} $\mathcal{M}$ \emph{iff }$(I$ \emph{is a
strong cut of }$\mathcal{M}$ \emph{and} $I\prec _{\Sigma _{1}}\mathcal{M)}$
\textit{iff\ }$I=\mathrm{Fix}(j)$ \emph{for some proper initial
self-embedding }$j$ \emph{of} $\mathcal{M}$.\textit{\smallskip }

\noindent \textbf{(c)}\textit{\ }$\mathrm{Fix}(j)=K^{1}(\mathcal{M})$ \emph{%
for some} \emph{self-embedding }$j$ \emph{of} $\mathcal{M}$ \emph{iff} $%
\mathbb{N}$ \textit{is a strong cut in }$\mathcal{M}$ \textit{iff} $\mathrm{%
Fix}(j)=K^{1}(\mathcal{M})$ \emph{for some proper initial self-embedding }$j$
\emph{of} $\mathcal{M}$.\textit{\medskip }

The plan of the paper is as follows: Section 2 reviews preliminaries;
Section 3 establishes some useful basic results about self-embeddings; and
Sections 4, 5, and 6 are respectively devoted to the proofs of parts (a),
(b), and (c) of Theorem 1.2. Some further results and open questions are
presented in Section 7.\medskip

\noindent \textbf{Acknowledgments.}~Saeideh Bahrami's research was partially
supported by the Iranian Ministry of Science, Research \& Technology, and
the Department of Philosophy, Linguistics \& Theory of Science of the
University of Gothenburg through funds which facilitated her three-month
visit to Gothenburg during 2016. Ali Enayat is indebted to Volodya Shavrukov
for playing a pivotal role in the inception of this paper since the
rudimentary forms of some of the results here were obtained by Enayat in the
course of brainstorming email discussions with Volodya during the winter and
spring months of 2012. These discussions also led to a number of questions,
which were eventually answered in this paper. Both authors are also grateful
to Costas Dimitracopoulos, Paul Gorbow, and especially Tin Lok Wong and the
anonymous referee for their assistance in weeding out infelicities in
earlier drafts of this paper.\textbf{\bigskip }

\begin{center}
\textbf{2.~PRELIMINARIES}\bigskip
\end{center}

In this section we review definitions, conventions, and known results that
will be utilized in this paper.

\begin{itemize}
\item The language of first order arithmetic, $\mathcal{L}_{A}$, is $%
\{+,\cdot ,\mathrm{S}(x),<,0\}.$ $\mathrm{PA}^{-}$ is the $\mathcal{L}_{A}$%
-theory describing the non-negative parts of discrete ordered rings as in
\cite{Kaye's text}. For a language $\mathcal{L}\supseteq \mathcal{L}_{A},$ $%
\mathrm{PA}(\mathcal{L})$ is $\mathrm{PA}^{-}$ augmented by the induction
scheme for all $\mathcal{L}$-formulae. We write $\mathrm{PA}$ for $\mathrm{PA%
}(\mathcal{L}_{A})$; when $\mathcal{L}$ is clear from the context, we shall
follow a common practice from the literature and use $\mathrm{PA}^{\ast }$
to refer to $\mathrm{PA}(\mathcal{L})$.

\item $M$, $M^{\ast },$ $M_{0}$, etc.~denote (respectively) the universes of
discourse of structures $\mathcal{M}$, $\mathcal{M}^{\ast },$ $\mathcal{M}%
_{0},$ etc. Given an $\mathcal{L}$-structure $\mathcal{M}$ and a class $%
\Gamma $ of $\mathcal{L}$-formulae, $\mathrm{Th}_{\Gamma }(\mathcal{M})$ is
the collection of sentences in $\Gamma $ that hold in $\mathcal{M}$. Also,
we write $\mathrm{Th}_{\exists }(\mathcal{M})$ for the collection of
existential sentences that hold in $\mathcal{M}$ (an existential formula is
of the form $\exists x_{0}\cdot \cdot \cdot \exists x_{k-1}\ \varphi $ for
quantifier-free $\varphi ).$

\item The meta-theoretic set of natural numbers is here denoted by $\omega $%
, and we use the notation $\left( a_{i}:i<s\right) $, where $s\in \omega $
or $s=\omega ,$ to refer to meta-theoretic sequences of finite or infinite
length. Given a model $\mathcal{M}$ of $\mathrm{PA}^{-}$, $\mathbb{N}$ is
the initial segment consisting of the standard elements of $\mathcal{M}$.
Also, given $s,i,$ and $a$ in $\mathcal{M}$, we write $\left( s\right)
_{i}=a $ to express the fact that $a$ is the $i$-th member of the sequence
canonically coded by $s$ in $\mathcal{M}$. In this context, we write $%
\left\langle a_{i}:i<r\right\rangle $ to refer to the object $s$ in $%
\mathcal{M}$ such that $s$ is the canonical code in $\mathcal{M}$ of a
sequence of length $r$ such that $\left( s\right) _{i}=a_{i}$ for each $i<r$%
. It is well-known \cite[Prop.~1.4.1]{Kossak-Schmerl} that we can arrange a
canonical coding such that if $s=\left\langle a_{i}:i<r\right\rangle $ and $%
a_{i}<b$ for all $i<r$, then $s\leq 2^{\left( r+b+1\right) ^{2}}$.

\item For a language $\mathcal{L}\supseteq \mathcal{L}_{A},$ $\Sigma _{0}(%
\mathcal{L})=\Pi _{0}(\mathcal{L})=\Delta _{0}(\mathcal{L})=$ the class of $%
\mathcal{L}$-formulae all of whose quantifiers are of the form $\exists x<t$
$\varphi $ or $\forall x<t\ \varphi $, where $t$ is an $\mathcal{L}$-term; $%
\Sigma _{n+1}(\mathcal{L})$ consists of formulae of the form $\exists
x_{0}\cdot \cdot \cdot \exists x_{k-1}\ \varphi $, where $\varphi \in \Pi
_{n}(\mathcal{L});$ and $\Pi _{n+1}(\mathcal{L})$ consists of formulae of
the form $\forall x_{0}\cdot \cdot \cdot \forall x_{k-1}\ \varphi $, where $%
\varphi \in \Sigma _{n}(\mathcal{L}).$ Here $k$ ranges over $\omega $, with
the understanding that $k=0$ corresponds to an empty block of quantifiers.
When $\mathcal{L}=\mathcal{L}_{A}$ we write $\Sigma _{n}$ and $\Pi _{n}$ for
$\Sigma _{n}(\mathcal{L})$ and $\Pi _{n}(\mathcal{L})$ (respectively).

\item For $n\in \omega ,\mathrm{I}\Sigma _{n}(\mathcal{L})$ is the fragment
of $\mathrm{PA}$ with the induction scheme limited to $\Sigma _{n}(\mathcal{L%
})$-formulae. The $\Sigma _{n}(\mathcal{L})$-Collection Scheme, denoted $%
\mathrm{B}\Sigma _{n}(\mathcal{L})$, consists of the universal closure of
formulae of the following form where $\varphi \in \Sigma _{n}(\mathcal{L})$
and $\varphi $ is allowed to have undisplayed parameters:
\end{itemize}

\begin{center}
$\left[ \forall x<v\ \exists y\ \varphi (x,y)\right] \rightarrow \exists z\ %
\left[ \forall x<v\ \exists y<z\ \varphi (x,y)\right] .$
\end{center}

\begin{itemize}
\item Given a theory $T$, and a class $\Gamma $ of formulae, $\Gamma ^{T}$
is the class of formulae that are $T$-provably equivalent to some formula in
$\Gamma $. It is well-known \cite[Ch.~7]{Kaye's text} that $\Sigma _{n}^{T}$
and $\Pi _{n}^{T}$ are both closed under bounded quantification,
disjunction, and conjunction for $T=\mathrm{I}\Delta _{0}+\mathrm{B}\Sigma
_{n}$.

\item For models $\mathcal{M}$ and $\mathcal{N}$ of $\mathcal{L}_{A}$, we
say that $\mathcal{N}$ \textit{end extends} {$\mathcal{M}$} (equivalently: $%
\mathcal{M}$ is an \textit{initial} submodel of $\mathcal{N}$), if $\mathcal{%
M}$ is a submodel of $\mathcal{N}$ and $a<b$ for every $a\in M,$ and $b\in
N\backslash M.$ For a class $\Gamma $ of $\mathcal{L}$-formulae we write $%
\mathcal{M}\prec _{\Gamma }\mathcal{N}$ if $\mathcal{N}$ is a $\Gamma $-%
\textit{elementary extension} of $\mathcal{M}$, i.e., $\Gamma $-formulae
with parameters in $\mathcal{M}$ are absolute in the passage between $%
\mathcal{M}$ and $\mathcal{N}$. An\textit{\ embedding }$\mathcal{M}$ into $%
\mathcal{N}$ is an isomorphism $j$ between $\mathcal{M}$ and a submodel of $%
\mathcal{N}$; such an embedding $j$ is said to be an \textit{initial
embedding} if the range of $j$ is an initial segment of $\mathcal{N}$. An
\emph{initial self-embedding} of $\mathcal{M}$ is an initial embedding of $%
\mathcal{M}$ into itself. A self-embedding $j$ is \textit{proper} if $j$ is
not surjective (equivalently, if $j$ is not an automorphism), otherwise $j$
is said to be \textit{improper}. Also, we say that a self-embedding $j$ is
\textit{trivial} if $j$ is the identity map on $\mathcal{M}$; otherwise $j$
is \textit{nontrivial}. Under these definitions, every automorphism of $%
\mathcal{M}$ is an improper initial self-embedding; and every proper
self-embedding is nontrivial.

\item $\mathrm{ACA}_{0}$ is the well-known subsystem of second order
arithmetic with the comprehension scheme limited to formulae with no second
order quantifiers, as in \cite{Steve-book}. Models of $\mathrm{ACA}_{0}$ are
of the two-sorted form $(\mathcal{M},\mathcal{A})$, where $\mathcal{A}$ is a
family of subsets of $M$, $(\mathcal{M},S)_{S\in \mathcal{A}}\models \mathrm{%
PA}^{\ast }$, and $\mathcal{A}$\ is closed under arithmetical definability. $%
\mathrm{WKL}_{0}$ is a subsystem of $\mathrm{ACA}_{0}$ whose models are of
the form $(\mathcal{M},\mathcal{A}),$ where $(\mathcal{M},\mathcal{A})$
satisfies (1) Induction for $\Sigma _{1}^{0}$ formulae (where $\Sigma
_{1}^{0}$ is the family of $\Sigma _{1}(\mathcal{L(A)})$ formulae with no
second order quantifier); (2) Comprehension for $\Delta _{1}^{0}$-formulae;
and (3) Weak K\"{o}nig's Lemma (which asserts that every infinite subtree of
the full binary tree has an infinite branch).
\end{itemize}

The following result is due to Paris and Pudl\'{a}k; it refines Bennett's
celebrated result stating that the graph of the exponential function $%
y=2^{x} $ is definable by a $\Delta _{0}$-predicate in the standard model of
arithmetic. See \cite[Sec.~V3(c)]{Hajek-Pudlak text} for further
detail.\medskip

\noindent \textbf{2.1.}~\textbf{Theorem}.~(Paris, Pudl\'{a}k) \emph{There is
a}\textbf{\ }$\Delta _{0}$-\emph{formula }$\mathrm{Exp}(x,y)$ \emph{such that%
} $\mathrm{I}\Delta _{0}$\textit{\ proves the following three statements}:%
\textit{\smallskip }

\noindent \textbf{(a)} $\forall x\exists ^{\leq 1}y$ $\mathrm{Exp}(x,y)$.%
\textit{\smallskip }

\noindent \textbf{(b)} $\forall x(\exists y$ $\mathrm{Exp}(x,y)\rightarrow
\forall z<x\ \exists y$ $\mathrm{Exp}(z,y))$.\textit{\smallskip }

\noindent \textbf{(c)} $\forall x\forall y$ $(\mathrm{Exp}(x,y)\rightarrow $
$\mathrm{Exp}(x+1,2y))$.\medskip

\begin{itemize}
\item $\mathrm{I}\Delta _{0}+\mathrm{Exp}$ is the extension of $\mathrm{I}%
\Delta _{0}$ obtained by adding the axiom $\mathrm{Exp}$, where $\mathrm{Exp}%
:=\forall x\exists y\ \mathrm{Exp}(x,y).$ The theory $\mathrm{I}\Delta _{0}+%
\mathrm{Exp}$ might not appear to be particularly strong since it cannot
even prove the totality of the superexponential function, but experience has
shown that it is a remarkably robust theory that is able to prove an
extensive array of theorems of number theory and finite combinatorics.

\item A \textit{cut} $I$ of a model $\mathcal{M}$\ of $\mathrm{PA}^{-}$ is
an initial segment of $\mathcal{M}$ with no last element. We write $m<I$,
where $m\in M$, to indicate that some member of $I$\ exceeds $m$. Similarly,
we write $I<m$ to indicate that every member of $I$ is below $m$. When a cut
$I$\ is closed under multiplication (and therefore under addition as well),
we shall use $I$ also to refer to the submodel of the ambient model whose
universe is $I$.
\end{itemize}

The following result is folklore; the verification that $\mathrm{I}\Delta
_{0}$ holds in $I$ is done by a routine induction on the length of $\Delta
_{0}$-formulae; see \cite[Prop.~10.5 ($n=1$)]{Kaye's text} for a proof that $%
\mathrm{B}\Sigma _{1}$ holds in $I$.\medskip

\noindent \textbf{2.2.}~\textbf{Theorem.}~\textit{If }$I$\textit{\ is a
proper cut of a model of }$\mathrm{I}\Delta _{0}$\textit{\ and} $I$ \textit{%
is closed under multiplication, then }$I\models \mathrm{I}\Delta _{0}+%
\mathrm{B}\Sigma _{1}.$ \medskip

\begin{itemize}
\item We will use $E$ to denote \textit{Ackermann's membership relation}
defined by: $xEy$ iff the $x$-th bit of the binary expansion of $y$ is a 1.
It is well-known that within $\mathrm{I}\Delta _{0}+\mathrm{Exp}$ the
formula $xEy$ is equivalent to a $\Delta _{0}$-formula. A subset $X$\ of $M$
is coded in $\mathcal{M}$ iff for some $m\in M$,
\end{itemize}

\begin{center}
$X=\left( m_{E}\right) ^{\mathcal{M}}:=\{x\in M:\mathcal{M}\models xEm\}.$
\end{center}

\begin{itemize}
\item Given $m\in M,$ $\underline{m}^{\mathcal{M}}:=\{x\in M:x<^{\mathcal{M}%
}m\}$. Note that $\underline{m}$ is coded in $\mathcal{M}\models \mathrm{I}%
\Delta _{0}$ provided $2^{m}$ exists in $\mathcal{M}$. When $\mathcal{M}$ is
clear from the context, we simply write $\underline{m}$ for $\underline{m}^{%
\mathcal{M}}.$

\item $X$ is \textit{piece-wise coded} in $\mathcal{M}$ if $\underline{m}%
\cap X$ is coded in $\mathcal{M}$ for each $m$ in $\mathcal{M}$.

\item For a cut $I$ of $\mathcal{M}$, $\mathrm{SSy}_{I}(\mathcal{M})$ is the
family consisting of sets of the form $S\cap I$, where $S$ is a subset of $M$
that is coded in $\mathcal{M}$, i.e.,
\begin{equation*}
\mathrm{SSy}_{I}(\mathcal{M})=\{\left( c_{E}\right) ^{\mathcal{M}}\cap
I:c\in M\}.
\end{equation*}%
When $I=\mathbb{N}$, we shall write the commonly used notation $\mathrm{SSy}(%
\mathcal{M})$ instead of $\mathrm{SSy}_{\mathbb{N}}(\mathcal{M}).$ It is
well-known \cite[Cor.~3.1]{Paola+Julia} that $\left( \mathbb{N},\mathrm{SSy}(%
\mathcal{M})\right) \models \mathrm{WKL}_{0}$ for a nonstandard $\mathcal{M}%
\models \mathrm{I}\Delta _{0}$; in particular $\mathrm{SSy}(\mathcal{M})$ is
a Boolean algebra and closed under Turing reducibility.

\item $\Delta _{0}(\Sigma _{n})$ is the class of $\mathcal{L}_{A}$-formulae
obtained by closing the class of $\Sigma _{n}$-formulae under Boolean
connectives and bounded quantifiers.

\item For a formula $\varphi (x_{1},\cdot \cdot \cdot ,x_{k})$ whose free
variables are ordered as shown, we write $\varphi ^{\mathcal{M}}$ for $%
\left\{ \left( m_{1},\cdot \cdot \cdot ,m_{k}\right) \in M^{k}:\mathcal{%
M\models \varphi }\left( m_{1},\cdot \cdot \cdot ,m_{k}\right) \right\} $.

\item Given a class $\Gamma $ of formulae, the $\Gamma $-Strong Collection
Scheme, here denoted $\mathrm{B}^{+}\Gamma $, consists of the universal
closure of formulae of the following form, where $\varphi (x,y)\in \Gamma $
and $\varphi $ is allowed to have undisplayed parameters:
\end{itemize}

\begin{center}
$\exists z\forall x<v\left[ \exists y\ \varphi (x,y)\rightarrow \exists y<z\
\varphi (x,y)\right] .$
\end{center}

\begin{itemize}
\item $\mathrm{Sat}_{_{\Sigma _{n}}}$ is the $\mathcal{L}_{\mathrm{A}}$%
-formula defining the satisfaction predicate for $\Sigma _{n}$-formulae for
an ambient model satisfying $\mathrm{I}\Delta _{0}+\mathrm{Exp}$. It is
well-known that $\mathrm{Sat}_{_{\Sigma _{n}}}\in \Sigma _{n}^{\mathrm{I}%
\Sigma _{1}}$ for each positive $n\in \omega ,$ and $\mathrm{Sat}_{_{\Sigma
_{0}}}\in \Sigma _{1}^{\mathrm{I}\Sigma _{1}}$ \cite[Thm.~1.75]{Hajek-Pudlak
text}.
\end{itemize}

The following theorem collects together a number of important properties of
models of $\mathcal{M}\models \mathrm{I}\Sigma _{n}$; see \cite[Ch.~I]%
{Hajek-Pudlak text} for an exposition.\smallskip

\noindent \textbf{2.3.}~\textbf{Theorem.}$~$\textit{If}$\mathcal{\ }n\in
\mathbb{\omega }$, $\mathcal{M}\models \mathrm{I}\Sigma _{n},$ \textit{and} $%
\varphi $\textit{\ is a unary }$\Delta _{0}(\Sigma _{n})$\textit{-formula }$%
\varphi (x,a)$\textit{,} \textit{where }$a$\textit{\ is a parameter from }$%
\mathcal{M},$ \textit{then}:\smallskip

\noindent \textbf{(a) }\textit{If} $n>0$, \textit{then} $\mathcal{M}\models
\mathrm{B}^{+}(\Sigma _{n})$.\smallskip

\noindent \textbf{(b)} $\varphi ^{\mathcal{M}}$ \emph{is piece-wise coded in}
$\mathcal{M}$ \textit{if} $n>0$, \textit{or} \textit{if} $n=0$ \textit{and} $%
\mathcal{M}\models \mathrm{Exp}.$\smallskip

\noindent \textbf{(c) }[$\Delta _{0}(\Sigma _{n})$-Min] \textit{If }$\varphi
^{\mathcal{M}}$ \textit{is nonempty}, \textit{then }$\varphi ^{\mathcal{M}}$%
\textit{\ has a minimum element.}\smallskip

\noindent \textbf{(d)} [$\Delta _{0}(\Sigma _{n})$-Max] \textit{If }$\varphi
^{\mathcal{M}}$\textit{\ is nonempty and bounded in }$\mathcal{M}$\textit{,
then }$\varphi ^{\mathcal{M}}$\textit{\ has a maximum element.}\smallskip

\noindent \textbf{(e)} [$\Delta _{0}(\Sigma _{n})$-Overspill] \textit{If }$%
\varphi ^{\mathcal{M}}$\textit{\ includes a proper cut }$I$ \textit{of} $%
\mathcal{M}$\textit{, then} $\underline{m}\subseteq \varphi ^{\mathcal{M}}$
\textit{for some }$m>\ I.$\smallskip

\noindent \textbf{(f) }[$\Delta _{0}(\Sigma _{n})$-PHP] \textit{If} $n>0$
\textit{and} $\varphi ^{\mathcal{M}}$\textit{\ is the graph of a function }$%
f $\textit{\ from }$\underline{m+1}$\textit{\ into }$\underline{m}$\textit{,
then }$f$ \textit{is not one-to-one.}\smallskip

\noindent \textbf{2.3.1.}$~$\textbf{Remark.}$~$Suppose $\mathcal{M}$ is a
nonstandard model of $\mathrm{I}\Sigma _{n}$ for $n>0$, and $p(x)$ is a
collection of formulae $\varphi (x,a)$ (where $a$ is a parameter in $%
\mathcal{M}$) such that (1) $p(x)$ is a $\Sigma _{n}$-type (i.e., every $%
\varphi \in p(x)$ is a $\Sigma _{n}$-formula); or (2) $p(x)$ is a short $\Pi
_{n}$-type (i.e., $p(x)$ includes the formula $x<\left( a\right) _{i}$ for
some $i\in \omega ,$ and every $\varphi \in p(x)$ is a $\Pi _{n}$-formula).
Then using part (e) of Theorem 2.3 (with $I=\mathbb{N}$), and the fact that $%
\mathrm{Sat}_{_{\Sigma _{n}}}$ has a $\Sigma _{n}$-description in $\mathcal{M%
}$ it is routine to verify that if $p(x)$ is coded in $\mathcal{M}$ (i.e., $%
\left\{ \ulcorner \varphi (x,y)\urcorner :\varphi \in p(x)\right\} \in
\mathrm{SSy}(\mathcal{M})$) and $p(x)$ is finitely realizable in $\mathcal{M}
$, then $p(x)$ is realized in $\mathcal{M}$.

\begin{itemize}
\item Given a class $\Gamma $ of formulae and $\mathcal{M}\models \mathrm{PA}%
^{-},$ $m\in M$ is said to be $\Gamma $-\textit{definable} in $\mathcal{M}$
if $\{m\}=\gamma ^{\mathcal{M}}$ for some unary $\gamma (x)\in \Gamma ;$ and
$m$ is $\Gamma $-\textit{minimal} in $\mathcal{M}$ if there is unary $\gamma
(x)\in \Gamma $ such that $m$ is the first element of $\gamma ^{\mathcal{M}%
}. $ Note that $m$ is $\Delta _{0}$-definable iff $m$ is $\Delta _{0}$%
-minimal. In general, if $m$ is $\Gamma $-definable then $m$ is $\Gamma $%
-minimal (but not conversely).

\item Given $\mathcal{M}\models \mathrm{PA}^{-},$ $K^{n}(\mathcal{M})$ is
the submodel of $M$ whose universe consists of all $\Sigma _{n}$-definable
elements of $\mathcal{M}$. The following result was originally proved by
Paris \& Kirby \cite[Prop.~8]{Jeff and Laurie}; see \cite[Ch.~IV]%
{Hajek-Pudlak text} for an expository account.
\end{itemize}

\noindent \textbf{2.4.}$~$\textbf{Theorem.}$~$(Paris \& Kirby) \emph{Suppose
}$n\in \mathbb{\omega }$ \emph{and} $\mathcal{M}\models \mathrm{I}\Sigma
_{n+1}.$\textit{\smallskip }

\noindent \textbf{(a) }$K^{n+1}(\mathcal{M})\prec _{\Sigma _{n+1}}\mathcal{M}
$.\textit{\smallskip }

\noindent \textbf{(b)} $K^{n+1}(\mathcal{M})\models \mathrm{I}\Sigma
_{n}+\lnot \mathrm{B}\Sigma _{n+1}$, \textit{if} $K^{n+1}(\mathcal{M})$
\textit{is nonstandard}$.$

\begin{itemize}
\item Given a cut $I$ of $\mathcal{M}$, $I$ is said to be a \textit{strong
cut} of $\mathcal{M}$ if, for each function $f$ whose graph is coded in $%
\mathcal{M}$ and whose domain includes $I,$\ there is some $s$\ in $M$ such
that for all $m\in I,$ $f(m)\notin I\ $iff $s<f(m).$ Paris \& Kirby proved
that strong cuts of models of $\mathrm{PA}$ are themselves models of $%
\mathrm{PA}$ \cite[Prop.~8]{Jeff and Laurie}. Indeed, their proof shows the
following more general result (see \cite[Sec.~7.3]{Kossak-Schmerl} or \cite[%
Lem.~A.4]{Me Tehran}).
\end{itemize}

\noindent \textbf{2.5.}$~$\textbf{Theorem.}$~$(Paris \& Kirby)\textit{\ The
following are equivalent for a proper cut }$I$\textit{\ of }$\mathcal{M}%
\models \mathrm{I}\Delta _{0}$:\textit{\smallskip }

\noindent \textbf{(1)} $I$ \textit{is a strong cut of} $\mathcal{M}$.\textit{%
\smallskip }

\noindent \textbf{(2)} $(I,\mathrm{SSy}_{I}(\mathcal{M}))\models \mathrm{ACA}%
_{0}.$\medskip

\begin{itemize}
\item Given a linearly ordered structure $\mathcal{K}$, let $\mathrm{Aut}(%
\mathcal{K})$ be the automorphism group of $\mathcal{K}$; $\mathrm{SE}(%
\mathcal{K})$ be the semi-group of self-embeddings of $\mathcal{K}$; $%
\mathrm{ISE}(\mathcal{K})$ be the semi-group of initial self-embeddings of $%
\mathcal{K}$, and $\mathrm{PISE}(\mathcal{K})$ be the semi-group of all
proper initial self-embeddings of $\mathcal{K}$ (all under composition).
Also, a self-embedding $j$ of $\mathcal{K}$ is \textit{contractive} iff $%
j(a)\leq a$ for all $a\in K.$
\end{itemize}

Theorem 2.6 below summarizes some remarkable results of Gaifman \cite[%
Thm.~4.9-4.11]{Gaifman}; his results were couched in terms of arbitrary
models of $\mathrm{PA}(\mathcal{L})$\ for countable $\mathcal{L}$ and are
proved using the technology of `minimal types'.\footnote{%
Note that if $(\mathcal{M},\mathcal{A})\models \mathrm{ACA}_{0}$, then the
expansion $(\mathcal{M},A)_{A\in \mathcal{A}}$ of $\mathcal{M}$ is a model
of $\mathrm{PA}(\mathcal{L})$, where $\mathcal{L}$ is the extension of $%
\mathcal{L}_{A}$ by predicate symbols for each $A\in \mathcal{A}.$ Moreover,
the collection of subsets of $M$ that are parametrically definable in $(%
\mathcal{M},A)_{A\in \mathcal{A}}$ coincides with $\mathcal{A}$.} A
streamlined proof of part (a) and the right-to-left direction of part (e)
appears in \cite[Thm.~B]{Me Tehran}. Part (h) of Theorem 2.6 seems to be
absent in Gaifman's paper; but a proof can be found in \cite[Thm.~3.3.8(c)]%
{Me unified}; the proof there is written for $j\in \mathrm{Aut}(\mathbb{L})$%
, but the reasoning carries over for $j\in \mathrm{SE}(\mathbb{L})$.\medskip

\noindent \textbf{2.6.}$\mathrm{~}$\textbf{Theorem.}$~$(Gaifman) \textit{%
Suppose }$(\mathcal{M},\mathcal{A})$\textit{\ is a countable model of }$%
\mathrm{ACA}_{0}$.\textit{\ Given any linear order} $\mathbb{L}$, \textit{%
there is} $\mathcal{N}_{\mathbb{L}}\succ _{\mathrm{end}}\mathcal{M}$ \textit{%
and an isomorphic copy} $\mathbb{L}^{\prime }=\{c_{l}:l\in \mathbb{L\}}$
\textit{of} $\mathbb{L}$ \textit{in} $N_{\mathbb{L}}\backslash M$, \textit{%
along with a composition preserving embedding }$j\mapsto \widehat{j}$ of $%
\mathrm{SE}(\mathbb{L})$\ \textit{into} $\mathrm{SE}(\mathcal{N}_{\mathbb{L}%
} $) \textit{such that}:\textit{\smallskip }

\noindent \textbf{(a) }$\mathrm{SSy}_{M}(\mathcal{N}_{\mathbb{L}})=\mathcal{A%
}$ \textit{and }$M\subseteq \mathrm{Fix}(\widehat{j})$ \textit{for each} $%
j\in \mathrm{SE}(\mathbb{L})$; \textit{moreover} $M=\mathrm{Fix}(\widehat{j}%
) $ \textit{iff} $j$ \textit{is fixed point free}.\textit{\smallskip }

\noindent \textbf{(b) }$\widehat{j}$\textit{\ is an \textbf{elementary}
self-embedding} \textit{of} $\mathcal{N}_{\mathbb{L}}$ \textit{for each} $%
j\in \mathrm{SE}(\mathbb{L})$.\textit{\smallskip }

\noindent \textbf{(c) }$\mathbb{L}^{\prime }$ \textit{is downward cofinal in
}$N_{\mathbb{L}}\backslash M$ \textit{if} $\mathbb{L}$ \textit{has no first
element. \smallskip }

\noindent \textbf{(d)} \textit{For any} $l_{0}\in \mathbb{L}$, $l_{0}$
\textit{is a strict upper bound for} $j(\mathbb{L})$ \textit{iff} $c_{l_{0}}$
\textit{is a strict upper bound for} $\widehat{j}(N_{\mathbb{L}})$.\smallskip

\noindent \textbf{(e) }$\widehat{j}\in \mathrm{Aut}(\mathcal{N}_{\mathbb{L}%
}) $ \textit{iff} $j\in $ $\mathrm{Aut}(\mathbb{L)}$.\textit{\smallskip }

\noindent \textbf{(f) }$\widehat{j}\in \mathrm{ISE}(\mathcal{N}_{\mathbb{L}%
}) $ \textit{iff} $j\in $ $\mathrm{ISE}(\mathbb{L)}.$\textit{\smallskip }

\noindent \textbf{(g) }$\widehat{j}\in \mathrm{PISE}(\mathcal{N}_{\mathbb{L}%
})$ \textit{iff} $j\in $ $\mathrm{PISE}(\mathbb{L)}$.\textit{\smallskip }

\noindent \textbf{(h) }$\widehat{j}$ \textit{is contractive iff} $j$ \textit{%
is contractive}.\textit{\medskip }

The following is Smory\'{n}ski's refinement of Friedman's embedding theorem.
The proof is outlined in \cite[Thm.~3.9]{Smorynski-Utrecht}, and given in
detail in \cite[Thm.~2.4]{Smorynski-JSL} (Smory\'{n}ski proved his result
for countable nonstandard models of $\mathrm{PA}$; but the proof readily
goes through for countable nonstandard models of I$\Sigma _{1})$.$\medskip $

\noindent \textbf{2.7.}$~$\textbf{Theorem.}$~$(Smory\'{n}ski) \textit{Suppose%
} $\mathcal{M}$ \textit{and} $\mathcal{N}$ \textit{are countable nonstandard
models of} $\mathrm{I}\Sigma _{1}$. \textit{The following are equivalent}:%
{\Large \smallskip }

\noindent \textbf{(1)} \textit{There is an embedding of} $\mathcal{M}$
\textit{into} $\mathcal{N}.${\Large \smallskip }

\noindent \textbf{(2)} $\mathrm{SSy}(\mathcal{M})\subseteq \mathrm{SSy}(%
\mathcal{N})$ \textit{and} $\mathrm{Th}_{\Sigma _{1}}(\mathcal{M})\subseteq
\mathrm{Th}_{\Sigma _{1}}(\mathcal{N}).${\Large \smallskip }

\noindent \textbf{(3) }\textit{There is an embedding }$j$\textit{\ of }$%
\mathcal{M}$ \textit{into} $\mathcal{N}$ \textit{such that }$j(\mathcal{M})$
\textit{is a `mixed' submodel of} $\mathcal{N}$, \textit{i.e.,} $j(\mathcal{M%
})$ \textit{is neither cofinal in }$\mathcal{M}$ \textit{nor an initial
segment of} $\mathcal{N}$.{\Large \smallskip }

\noindent \textbf{2.7.1.}$~$\textbf{Remark.}$~$As noted by Smory\'{n}ski
\cite[p.~21]{Smorynski-Lectures (logic colloquium)} the condition $\mathrm{Th%
}_{\Sigma _{1}}(\mathcal{M})\subseteq \mathrm{Th}_{\Sigma _{1}}(\mathcal{N})$
in (2) above can be weakened to $\mathrm{Th}_{\exists }(\mathcal{M}%
)\subseteq \mathrm{Th}_{\exists }(\mathcal{N}),$ thanks to the MRDP Theorem.
The MRDP Theorem (due to Matijasevi\v{c}, Robinson, Davis, and Putnam)
states that every recursively enumerable set is Diophantine. As shown by
Dimitracopoulos and Gaifman \cite{Costas-Haim} the MRDP Theorem is provable
in $\mathrm{I\Delta }_{0}+\mathrm{Exp.}$\medskip

The next result is due to Wilkie (according to \cite{Smorynski-JSL}, where
it first appeared in print). Wilkie's result was formulated for countable
nonstandard models of $\mathrm{PA}$, but an inspection of the proof
presented in \cite[Thm.~12.6]{Kaye's text} makes it clear that the result
holds for countable nonstandard models of $\mathrm{I}\Sigma _{2}$.\medskip

\noindent \textbf{2.8.}$~$\textbf{Theorem.}$~$(Wilkie) \textit{Suppose} $%
\mathcal{M}$ \textit{and} $\mathcal{N}$ \textit{are countable nonstandard
models of} $\mathrm{I}\Sigma _{2}$. \textit{The following are equivalent}:%
{\Large \smallskip }

\noindent \textbf{(1)} \textit{For each} $a\in N$ \textit{there is a proper
initial embedding} $j$ \textit{of} $\mathcal{M}$ \textit{into} $\mathcal{N}$
\textit{such that} $a\in j(M).${\Large \smallskip }

\noindent \textbf{(2)} $\mathrm{SSy}(\mathcal{M})=\mathrm{SSy}(\mathcal{N})$
\textit{and} $\mathrm{Th}_{\Pi _{2}}(\mathcal{M})\subseteq \mathrm{Th}_{\Pi
_{2}}(\mathcal{N}).$\medskip

The following result of Ressayre \cite{Ressayre} shows that all countable
nonstandard models of $\mathrm{I}\Sigma _{1}$ carry proper initial
self-embeddings that pointwise fix any prescribed topped initial segment;
and $\mathrm{I}\Sigma _{1}$ is the weakest extension of $\mathrm{I}\Delta
_{0}$ with this property. The $(1)\Rightarrow (2)$ direction of Ressayre's
theorem is refined in Corollary 3.3.1 and Theorem 4.1; see Remarks 3.3.2 and
4.1.2 for more detail.\medskip

\noindent \textbf{2.9.}$\mathrm{~}$\textbf{Theorem.}$~$(Ressayre) \textit{%
The following are equivalent for a countable nonstandard }$\mathcal{M}%
\models \mathrm{I}\Delta _{0}$:\textit{\smallskip }

\noindent \textbf{(1)} $\mathcal{M}$ $\models \mathrm{I}\Sigma _{1}.$\textit{%
\smallskip }

\noindent \textbf{(2)} \textit{For each }$a\in M$\textit{, there is a proper
initial self-embedding }$j$\textit{\ of }$\mathcal{M}$\textit{\ such that }$%
j(m)=m$\textit{\ for each }$m\leq a.$\textbf{\bigskip }

\begin{center}
\textbf{3.~BASIC RESULTS}\textit{\bigskip }
\end{center}

In this section we establish a number of basic results about
self-embeddings. These results will also be useful in subsequent sections.
\medskip

\noindent \textbf{3.1.}$\mathrm{~}$\textbf{Theorem.}$\mathrm{~}$\emph{Suppose%
} $j$ \emph{is a self-embedding of }$\mathcal{M}\models \mathrm{I\Delta }%
_{0}+\mathrm{Exp}.$ \emph{Then} $K^{1}(\mathcal{M})\preceq _{\Sigma _{1}}%
\mathrm{Fix}(j)\preceq _{\Sigma _{1}}\mathcal{M}.\medskip $

Before presenting the proof of Theorem 3.1, we will establish two useful
lemmas.$\medskip $

\noindent \textbf{3.1.1.}$\mathrm{~}$\textbf{Lemma.}$~$\emph{If} $\mathcal{M}
$ \emph{and} $\mathcal{N}$ \emph{are both models of} $\mathrm{I\Delta }_{0}+%
\mathrm{Exp},$ \emph{and} $j$ \emph{is an embedding of }$\mathcal{M}$ \emph{%
into} $\mathcal{N}$, \emph{then} $j(\mathcal{M})\preceq _{\Delta _{0}}%
\mathcal{N}.$\medskip

\noindent \textbf{Proof.}$\mathrm{~}$If $j$ is an initial embedding, then
this follows from the basic fact that\textbf{\ }every submodel of $\mathcal{N%
}$ whose universe is a cut of $\mathcal{N}$ that is closed under
multiplication (and therefore addition) is a $\Delta _{0}$-elementary
submodel of $\mathcal{N}$. For the general case, this follows from the
provability of the MRDP Theorem in models of $\mathrm{I\Delta }_{0}+\mathrm{%
Exp,}$ since if $\mathcal{N}_{0}$ is a submodel of $\mathcal{N}$, where both
$\mathcal{N}_{0}$ and $\mathcal{N}$ are models of $\mathrm{I\Delta }_{0}+%
\mathrm{MRDP}$, then $\mathcal{N}_{0}\preceq _{\Delta _{0}}\mathcal{N}$.
\hfill $\square \medskip $

\noindent \textbf{3.1.2.}$~$\textbf{Lemma.}$~$\emph{Suppose} $\mathcal{M}%
\models \mathrm{I\Delta }_{0}.\smallskip $

\noindent \textbf{(a) }\emph{If }$D$ \emph{is a} \emph{nonempty} $\Sigma
_{1} $-\emph{definable subset of} $\mathcal{M}$, \emph{then} \emph{there is
some }$d\in D$ \emph{such that} $d$ \emph{is }$\Delta _{0}$-\emph{minimal in}
$\left( \mathcal{M},m\right) $ \emph{for some} $\Delta _{0}$-\emph{minimal
element }$m$\emph{\ of }$\mathcal{M}$ \emph{with} $d<m$.$\mathtt{\smallskip }
$

\noindent \textbf{(b)} \emph{If} $d\in K^{1}(\mathcal{M})$, \emph{then} $d$
\emph{is} $\Delta _{0}$-\emph{minimal in} $\left( \mathcal{M},m\right) $
\emph{for some} $\Delta _{0}$-\emph{minimal element }$m$\emph{\ of }$%
\mathcal{M}$ \emph{with} $d<m.\mathtt{\smallskip }$

\noindent \textbf{(c) }\emph{If in addition }$\mathcal{M}\models \mathrm{Exp}
$, \emph{and} $j$ \emph{is a self-embedding of }$\mathcal{M}$ \emph{such that%
} $j(m)=m$, \emph{and} $d$ \emph{is} $\Delta _{0}$-\emph{minimal in} $(%
\mathcal{M},m)$, \emph{then} $j(d)=d.$ $\medskip $

\noindent \textbf{Proof.}$~$\textbf{(a) }Easy; suppose $D$ is definable by
the formula $\exists z\ \delta (x,z),$ where $\delta $ is $\Delta _{0}$. Let
$m$ be the first element in $\mathcal{M}$ such that $\delta (x,z)$ holds for
some $x$ and $z$ below $m$, and then let $d$ be the first element below $m$
such that $\delta (d,z)$ holds for some $z<m$. \medskip

\noindent \textbf{(b) }This follows immediately from part (a) by setting $%
D=\{d\}$.\medskip

\noindent \textbf{(c) }Suppose $\delta (x,y)$ is a $\Delta _{0}$-formula
such that:\smallskip

\noindent (1) $(\mathcal{M},m)\models d=\mu x$ $\delta (x,m)$,\smallskip

\noindent where $\mu $ is the least search operator. (1) coupled with the
assumption that $j$ is an isomorphism between $\mathcal{M}$ and $j(\mathcal{M%
})$ implies:\smallskip

\noindent (2) $\left( j(\mathcal{M)},j(m)\right) \models j(d)=\mu x$ $\delta
(x,j(m)).$ \smallskip

\noindent By putting (2) together with $j(\mathcal{M)\preceq }_{\Delta _{0}}%
\mathcal{M}$ (by Lemma 3.1.1) and the assumption $j(m)=m$ we have:\smallskip

\noindent (3) $(\mathcal{M},m)\models j(d)=\mu x$ $\delta (x,m)$.\smallskip

\noindent By putting (1) together with (3) we can now conclude that $j(d)=d$%
.\hfill $\square \medskip $

\noindent \textbf{Proof of Theorem 3.1.}$\mathrm{~}$Let us first establish $%
\mathrm{Fix}(j)\preceq _{\Sigma _{1}}\mathcal{M}$. By Tarski's test, it
suffices to show that for every $\Delta _{0}$-formula\emph{\ }$\delta (x,y)$%
, if $\mathcal{M}\models \exists x\ \delta (x,m)$ for some $m\in \mathrm{Fix}%
(j)$, then $\mathcal{M}\models \delta (d,m)$ for some $d\in \mathrm{Fix}(j).$
Let $D$ be defined in $\mathcal{M}$ as consisting of elements $x$ such that $%
\delta (x,m),$ and let $d$ be the least member of $D$. Then $d$ is $\Delta
_{0}$-minimal in $\left( \mathcal{M},m\right) $, and therefore $j(d)=d$ by
part (c) of Lemma 3.1.2. \medskip

To see that $K^{1}(\mathcal{M})\subseteq \mathrm{Fix}(j)$, suppose $d\in
K^{1}(\mathcal{M})$. Then by part (b) of Lemma 3.1.2 there is some $\Delta
_{0}$-minimal element $m$ of $\mathcal{M}$ such that $d$ is $\Delta _{0}$%
-minimal in $(\mathcal{M},m)$. Therefore by two applications of part (c) of
Lemma 3.1.2 we can obtain $j(m)=m$ and $j(d)=d$. Recall that $K^{1}(\mathcal{%
M})\preceq _{\Sigma _{1}}\mathcal{M}$ (by the $n=0$ case of Theorem 2.4),
and we have already verified that $K^{1}(\mathcal{M})\subseteq \mathrm{Fix}%
(j)$ and $\mathrm{Fix}(j)\preceq _{\Sigma _{1}}\mathcal{M}$. On the other
hand it can be easily seen that in general if $\mathcal{N}_{0}$ and $%
\mathcal{N}_{1}$ are $\Sigma _{1}$-elementary submodels of an $\mathcal{L}%
_{A}$-structure $\mathcal{N}$ with $N_{0}\subseteq N_{1}$, then $\mathcal{N}%
_{0}\preceq _{\Sigma _{1}}\mathcal{N}_{1}$. This completes the proof of $%
K^{1}(\mathcal{M})\preceq _{\Sigma _{1}}\mathrm{Fix}(j)$.\hfill $\square
\medskip $

\noindent \textbf{3.1.3.}$~$\textbf{Remark.}$~$It is easy to see, using part
(b) of Lemma 3.1.2, that $K^{1}(\mathcal{M})=\Delta _{1}^{\mathcal{M}}$;
i.e., the elements of $K^{1}(\mathcal{M})$ are precisely those elements of $%
\mathcal{M}$\textit{\ }that\textit{\ }are both $\Sigma _{1}$-definable and $%
\Pi _{1}$-definable in\textit{\ }$\mathcal{M}$. This observation dates back
to Mijajlovi\'{c} \cite{Zarko}.\medskip

The following result generalizes the $(a)\Rightarrow (b)$ direction of \cite[%
Thm.~A]{Me FM}, which corresponds to Theorem 3.2 when $j$ is a nontrivial
automorphism of $\mathcal{M}$. \medskip

\noindent \textbf{3.2.}$\mathrm{~}$\textbf{Theorem.}$\mathrm{~}$\emph{If}%
\textrm{\ }$\mathcal{M}\models \mathrm{I\Delta }_{0}$ \emph{and} $j$ \emph{%
is a \textbf{nontrivial} self-embedding of} $\mathcal{M}$ \emph{such that} $%
j(\mathcal{M})\preceq _{\Delta _{0}}\mathcal{M}$, \emph{then} $\mathrm{I}_{%
\mathrm{fix}}(j)\models \mathrm{I}\Delta _{0}+\mathrm{B}\Sigma _{1}+\mathrm{%
Exp.}\medskip $

\noindent \textbf{Proof.}$~$We first verify that $\mathrm{I}_{\mathrm{fix}%
}(j)$ is closed under the operations of the ambient structure $\mathcal{M}$.
Suppose $x$ and $y$ are elements of $\mathrm{I}_{\mathrm{fix}}(j)$ with $%
x\leq y$ and, without loss of generality, assume that $x$ and $y$ are both
nonstandard elements. Since $x+y<xy\leq y^{2}$, it suffices to show that $%
y^{2}\in \mathrm{I}_{\mathrm{fix}}(j).$ Observe that $I\Delta _{0}$ can
prove that any number $z<y^{2}$ can be written as $z=qy+r,$ where both $q$
and $r$ are less than $y$ (since the division algorithm can be implemented
in $\mathrm{I}\Delta _{0}$). Therefore,

\begin{center}
$j(z)=j(qy+r)=j(q)j(y)+j(r)=qy+r=z.$
\end{center}

\noindent This shows that $\mathrm{I}_{\mathrm{fix}}(j)$ is closed under the
operations of $\mathcal{M}\mathfrak{.}$ It is also clear by the definition
of $\mathrm{I}_{\mathrm{fix}}(j)$ and the assumption that $j$ moves some
element of $\mathcal{M}$ that $\mathrm{I}_{\mathrm{fix}}(j)$ is a proper cut
of $\mathcal{M}$. Hence $\mathrm{I}_{\mathrm{fix}}(j)\models \mathrm{I}%
\Delta _{0}+\mathrm{B}\Sigma _{1}$ by Theorem 2.2. \medskip

It remains to show that $\mathrm{Exp}$ holds in $\mathrm{I}_{\mathrm{fix}%
}(j) $. First we will show:\smallskip

\noindent $(\ast )$ If $a\in \mathrm{I}_{\mathrm{fix}}(j)$ and $2^{a}$ is
defined in $\mathcal{M}$, then $2^{a}\in \mathrm{I}_{\mathrm{fix}}(j)$.
\smallskip

\noindent To establish $(\ast )$, suppose $\mathcal{M}\models b<2^{a}$. Then
$\mathcal{M}\models b=\sum\limits_{i<c}2^{s_{i}}$, with $c\leq a$ and $%
s_{0}<\cdot \cdot \cdot <s_{c-1}<a$. Therefore $j(c)=c$ and $j(s_{i})=s_{i}$
for each $i<c$, because $a\in \mathrm{I}_{\mathrm{fix}}(j).$ So we have some
element $b^{\prime }\in j(M)$ such that:

\begin{center}
$j(\mathcal{M})\models
j(b)=\sum\limits_{i<j(c)}2^{j(s_{i})}=\sum\limits_{i<c}2^{s_{i}}=b^{\prime
}. $
\end{center}

\noindent But $j(\mathcal{M})\prec _{\Delta _{0}}\mathcal{M}$ by assumption,
and therefore the $j(\mathcal{M})$-binary representation of each element of $%
j(\mathcal{M})$ coincides with the $\mathcal{M}$-binary representation of
the same element since for a sequence $s=\left\langle s_{i}:i<c\right\rangle
$ in $j(\mathcal{M})$, where $c$ might be nonstandard, the statement $%
x=\sum\limits_{i<c}2^{s_{i}}$ is well-known to be expressible in $j(\mathcal{%
M})$ by a $\Delta _{0}$-formula $\delta (x,s,p)$ (where $p$ is some
sufficiently large parameter). This makes it clear that $b^{\prime }=b$.
Therefore $j(b)=b$ for each $b<2^{a}$; which in turn implies that $2^{a}\in
\mathrm{I}_{\mathrm{fix}}(j)$. \medskip

In light of $(\ast )$, the proof that $\mathrm{Exp}$ holds in $\mathcal{M}$
will be complete once we demonstrate that for all $a\in \mathrm{I}_{\mathrm{%
fix}}(j)$, $2^{a}$ is defined in $\mathcal{M}$. Indeed, we will establish
the slightly stronger result $(\ast \ast )$ below: \smallskip

\noindent $(\ast \ast )$ $\mathrm{I}_{\mathrm{fix}}(j)\subsetneq J,$ where $%
J:=\{x\in M:\mathcal{M}\models \exists y(2^{x}=y)\}.$\smallskip

\noindent In order to verify $(\ast \ast )$, first let $P:=\{y\in M:\mathcal{%
M}\models \exists x(2^{x}=y)\}$, and note that: \smallskip

\noindent $(1)$ $P$ is unbounded in $\mathcal{M}$,\smallskip

\noindent since otherwise by putting the fact that the graph of the
exponential function is $\Delta _{0}$-definable in $\mathcal{M}$ (Theorem
2.1) together with the veracity of $\Delta _{0}$-Max in $\mathcal{M}$
(Theorem 2.3(d)), there would have to be a last power of 2 in $\mathcal{M}$,
which is impossible. Next, note that if $(\ast \ast )$\ fails,
then:\smallskip

\noindent $(2)$ $J\subseteq \mathrm{I}_{\mathrm{fix}}(j)$, \smallskip

\noindent because $J$ is an initial segment of $\mathcal{M}$ by Theorem
2.1(b). By putting (2) together with $(\ast )$ we obtain:\smallskip

\noindent $(3)$ $P\subseteq \mathrm{I}_{\mathrm{fix}}(j)$. \smallskip

\noindent But since $j$ is assumed to be nontrivial, there is some $c\in M$
such that $\mathrm{I}_{\mathrm{fix}}(j)<c$, and so by (3) $P$\ is bounded
above by $c$, which contradicts (1), and thereby concludes the proof of $%
(\ast \ast )$. \hfill $\square \medskip $

Theorem 3.3 below fine-tunes a result of H\'{a}jek \& Pudl\'{a}k \cite[%
Thm.~11]{Hajek-Pudlak paper}. Their result is the special case of Corollary
3.3.1 when $\mathcal{M}$ and $\mathcal{N}$, as well as the cut $I,$ are all
assumed to be models of $\mathrm{PA}$. \textit{\medskip }

\noindent \textbf{3.3.}$~$\textbf{Theorem.}$~$\emph{Suppose }$\mathcal{M}$
\emph{and} $\mathcal{N}$ \emph{are countable nonstandard models of} $\mathrm{%
I}\Sigma _{1}$ \emph{with} $c\in M$ \emph{and} $a,b\in N$. \emph{%
Furthermore, suppose} $I$ \emph{is a proper cut shared by} $\mathcal{M}$
\emph{and} $\mathcal{N}$ \emph{such that} $I$ \textit{is closed under
exponentiation. }\emph{The following are equivalent}:\textit{\smallskip }

\noindent \textbf{(1)} \emph{There is a proper initial embedding} $j:%
\mathcal{M}\rightarrow \mathcal{N}$ \emph{such that} $j(c)=a$, $j(M)<b$
\emph{and} $j(i)=i$ \emph{for all }$i\in I$.\textit{\smallskip }

\noindent \textbf{(2)} $\mathrm{SSy}_{I}(\mathcal{M})=\mathrm{SSy}_{I}(%
\mathcal{N})$, \emph{and for all}\textrm{\ }$i\in I$ \emph{and all }$\Delta
_{0}$-\emph{formulae} $\delta (x,y,z)$,\emph{\ if} $\mathcal{M}\models
\exists z\ \delta (i,c,z)$,\emph{\ then} $\mathcal{N}\models \exists z<b\
\delta (i,a,z).$\textit{\medskip }

\noindent \textbf{Proof.}$~(1)\Rightarrow (2)$ is easy and is left to the
reader so we will concentrate on $(2)\Rightarrow (1)\mathbf{.}$ Assume (2)
and fix an enumeration $\left( c_{k}:k<\omega \right) $ of $M;$ and an
enumeration $\left( d_{k}:k<\omega \right) $ of $N$ in which each element of
$N$ occurs infinitely often. The proof of (1) will be complete by setting $%
j(u_{k})=v_{k}$ once we have $\left( u_{k}:k<\omega \right) $ and $\left(
v_{k}:k<\omega \right) $ that satisfy the following four conditions:\medskip

\noindent (I) $M=\left\{ u_{k}:k<\omega \right\} .\smallskip $

\noindent (II) $\left\{ v_{k}:k<\omega \right\} $ is an initial segment of $%
\mathcal{N}$, and each $v_{k}<b.\smallskip $

\noindent (III) $u_{0}=a$ and $v_{0}=c.\smallskip $

\noindent (IV) For each positive $n<\omega $, each $i\in I$, and each $%
\Delta _{0}$-formula $\delta (x,\mathbf{y},z)$, where $\mathbf{y}=\left(
y_{r}:r<n\right) $, the following holds for $\mathbf{u}=\left(
u_{r}:r<n\right) $, and $\mathbf{v}=\left( v_{r}:r<n\right) $:

\begin{center}
$\mathcal{M}\models \exists z\ \delta (i,\mathbf{u},z)\Longrightarrow
\mathcal{N}\models \exists z<b\ \delta (i,\mathbf{v},z)$.
\end{center}

\noindent We will define finite tuples $\mathbf{u}_{m}=\left(
u_{r}:r<n_{m}\right) $ and , $\mathbf{v}_{m}=\left( v_{r}:r<n_{m}\right) $
from $M$ (and of the same length) by recursion on $m$ so that the following
condition is maintained through the recursion for all $m<\omega $:$%
\smallskip $

\noindent $(\ast _{m})\ $If $\mathcal{M}\models \exists z\ \delta (i,\mathbf{%
u}_{m},z)$, then $\mathcal{N}\models \exists z<b\ \delta (i,\mathbf{v}%
_{m},z) $, for all $i\in I,$ and each $\delta (x,\mathbf{y},z)\in \Delta
_{0} $, where $\mathbf{y}=\left( y_{r}:r<n_{m}\right) $.$\smallskip $

\noindent For $m=0$, we set $\mathbf{u}_{0}=\left( a\right) $ and $\mathbf{v}%
_{0}=\left( c\right) $, so $n_{0}=1.$ By (2) this choice of $\mathbf{u}_{0}$
and $\mathbf{v}_{0}$ satisfies $(\ast _{0})$. Let\textbf{\ }$\left\langle
\delta _{r}:r\in M\right\rangle $ be a canonical enumeration within $%
\mathcal{M}$ of all $\Delta _{0}$-formulae (e.g., as in \cite[Ch.~1]%
{Hajek-Pudlak text}). For $m\geq 0$, we may assume that there are $\mathbf{u}%
_{m}$ and $\mathbf{v}_{m}$ satisfying $(\ast _{m}).$ In order to construct $%
\mathbf{u}_{m+1}$ and $\mathbf{v}_{m+1}$ we distinguish between the case $%
m=2k$ (the $k$-th `forth' stage) and the case $m=2k+1$ (the $k$-th `back'
stage) as described below.\medskip

\noindent \texttt{CASE}\ $m=2k.$ In this case, if $c_{k}$ is already among
the elements listed in $\mathbf{u}_{m}$ we have nothing to do, i.e., in this
case $\mathbf{u}_{m+1}=\mathbf{u}_{m}$ and $\mathbf{v}_{m+1}=\mathbf{v}_{m}$%
. Otherwise, consider:

\begin{center}
$H:=\left\{ \left\langle r,i\right\rangle \in I:\mathcal{M}\models \exists
z\ \mathrm{Sat}_{\Delta _{0}}\left( \delta _{r}(i,\mathbf{u}%
_{m},c_{k},z)\right) \right\} $.
\end{center}

\noindent $H$ is the intersection of a $\Sigma _{1}$-definable subset of $%
\mathcal{M}$ with $I$, so $H\in \mathrm{SSy}_{I}(\mathcal{M})=\mathrm{SSy}%
_{I}(\mathcal{N})$. Therefore we can choose $h$ in $\mathcal{M}$ and $%
h^{\prime }$ in $\mathcal{N}$ such that:

\begin{center}
$H=I\cap \left( h_{E}\right) ^{\mathcal{M}}=I\cap \left( h_{E}^{\prime
}\right) ^{\mathcal{N}}.$
\end{center}

\noindent For each $p\in M$ and $q\in N$ define:

\begin{center}
$H_{p}:=\left( h_{E}\cap \underline{p}\right) ^{\mathcal{M}}$ and $%
H_{q}^{\prime }:=\left( h_{E}^{\prime }\cap \underline{q}\right) ^{\mathcal{N%
}}.$
\end{center}

\noindent Choose $h_{p}\in M$ and $h_{q}^{\prime }\in N$ such that $H_{p}$
is coded by $h_{p}$ in $\mathcal{M}$, and $H_{q}^{\prime }$ is coded by $%
h_{q}^{\prime }$ in $\mathcal{N}$. In light of the assumption that $I$ is
closed under exponentiation, we have:\smallskip

\noindent $(i)$ $h_{s}=h_{s}^{\prime }\in I$ for each $s\in I.$\smallskip

\noindent On the other hand, by definition:\smallskip

\noindent $(ii)$ $s\in I\Rightarrow \mathcal{M}\models \forall \left\langle
r,i\right\rangle \in h_{s}\ \exists z\ \mathrm{Sat}_{\Delta _{0}}(\delta
_{r}(i,\mathbf{u}_{m},c_{k},z))$.\smallskip

\noindent Putting $(ii)$ together with $\Sigma _{1}$-Collection in $\mathcal{%
M}$ yields: \smallskip

\noindent $(iii)$ $s\in I\Rightarrow \mathcal{M}\models \exists t\ \forall
\left\langle r,i\right\rangle \in h_{s}\ \exists z<t\ \mathrm{Sat}_{\Delta
_{0}}(\delta _{r}(i,\mathbf{u}_{m},c_{k},z))$.\smallskip

\noindent By quantifying out $c_{k}$ in $(iii)$ we obtain:\smallskip

\noindent $(iv)$ $s\in I\Rightarrow \mathcal{M}\models \overset{\varphi
(h_{s},\mathbf{u})}{\overbrace{\exists x\ \exists t\ \forall \left\langle
r,i\right\rangle \in h_{s}\ \exists z<t\ \mathrm{Sat}_{\Delta _{0}}(\delta
_{r}(i,\mathbf{u}_{m},x,z))}}$.\smallskip

\noindent Note that $\varphi (h_{s},\mathbf{u}_{m})$ can be written as a $%
\Sigma _{1}$-formula. Therefore by coupling our inductive hypothesis $(\ast
_{m})$ with $(i)$ and $(iii)$ we conclude:\smallskip

\noindent $(v)$ $s\in I\Rightarrow \ \mathcal{N}\models \exists x,t<b\
\forall \left\langle r,i\right\rangle \in h_{s}^{\prime }\ \exists z<t\
\mathrm{Sat}_{\Delta _{0}}(\delta _{r}(i,\mathbf{v}_{m},x,z)).$\smallskip

\noindent Finally, by $(v)$ together and $\Sigma _{1}$-Overspill in $%
\mathcal{N}$ there exists $p>I$ such that:\smallskip

\noindent $(vi)$ $\mathcal{N}\models \exists x,t<b\ \forall \left\langle
r,i\right\rangle \in h_{p}^{\prime }\ \exists z<t\ \mathrm{Sat}_{\Delta
_{0}}(\delta _{r}(i,\mathbf{v}_{m},x,z))$.\smallskip

\noindent Let $d$ be a witness in $\mathcal{N}$ to the $\exists x$ assertion
in $(vi)$, and let $\mathbf{u}_{m+1}=\left( \mathbf{u}_{m},c_{k}\right) $
and $\mathbf{v}_{m+1}=\left( \mathbf{v}_{n},d\right) .$ It is easy to see
using $(vi)$ that $(\ast _{m+1})$ holds with these choices of $\mathbf{u}%
_{m+1}$ and $\mathbf{v}_{m+1}.\medskip $

\noindent \texttt{CASE}\ $m=2k+1.$ If $d_{k}>\max \left( \mathbf{v}%
_{m}\right) $ we do nothing, i.e., we define $\mathbf{u}_{m+1}:=\mathbf{u}%
_{m}$ and $\mathbf{v}_{m+1}:=\mathbf{v}_{m}$. Otherwise, let:

\begin{center}
$L=\{\left\langle r,i\right\rangle \in I:\mathcal{N}\models \forall z\
\left( \mathrm{Sat}_{\Delta _{0}}(\delta _{r}(i,\mathbf{v}%
_{m},d_{k},z)\rightarrow b\leq z\right) \}$.
\end{center}

\noindent Since $L$ is the intersection of a $\Pi _{1}$-definable subset of $%
\mathcal{N}$ with $I$, $L\in \mathrm{SSy}_{I}(\mathcal{N})=\mathrm{SSy}_{I}(%
\mathcal{M})$. Therefore we can choose $l$ in $\mathcal{M}$ and $l^{\prime }$
in $\mathcal{N}$ such that:

\begin{center}
$L=I\cap \left( l_{E}\right) ^{\mathcal{M}}=I\cap \left( l_{E}^{\prime
}\right) ^{\mathcal{N}}.$
\end{center}

\noindent For each $p\in M$ and $q\in N$ define:

\begin{center}
$L_{p}:=\left( l_{E}\cap \underline{p}\right) ^{\mathcal{M}}$ and $%
L_{q}^{\prime }:=\left( l_{E}^{\prime }\cap \underline{q}\right) ^{\mathcal{N%
}}.$
\end{center}

\noindent Let $l_{p}\in M$ and $l_{q}^{\prime }\in N$ such that $L_{p}$ is
coded by $l_{p}$ in $\mathcal{M}$, and $L_{q}^{\prime }$ is coded by $%
l_{q}^{\prime }$ in $\mathcal{N}$. The closure of $I$ under exponentiation
makes it clear that:\smallskip

\noindent $(vii)$ $l_{s}=l_{s}^{\prime }\in I$ for each $s\in I.$\smallskip

\noindent We claim that for every $s\in I$ the following holds:\smallskip

\noindent $(viii)$ $\mathcal{M}\models \exists x\leq \mathrm{max}(\mathbf{u}%
_{m})\ \forall \left\langle r,i\right\rangle \in l_{s}\ \forall z\ \lnot
\mathrm{Sat}_{\Delta _{0}}(\delta _{r}(i,\mathbf{u}_{m},x,z))$.\smallskip

\noindent Suppose not, then for some $s\in I$:\smallskip

\noindent $(ix)$ $\mathcal{M}\models \forall x\leq \mathrm{max}(\mathbf{u}%
_{m})\ \exists \left\langle r,i\right\rangle \in l_{s}\ \exists z\ \mathrm{%
Sat}_{\Delta _{0}}(\delta _{r}(i,\mathbf{u}_{m},x,z)).$\smallskip

\noindent Thanks to $\Sigma _{1}$-collection in $\mathcal{M}$ and $(viii)$
we obtain:\smallskip

\noindent $(x)$ $\mathcal{M}\models \exists t\ \forall x\leq \mathrm{max}(%
\mathbf{u}_{m})\ \exists \left\langle r,i\right\rangle \in l_{s}\ \exists
z<t\ \mathrm{Sat}_{\Delta _{0}}(\delta _{r}(i,\mathbf{u}_{m},x,z)).$%
\smallskip

\noindent So by our inductive assumption $(\ast _{m})$ and $(x)$ we
have:\smallskip

\noindent $(xi)$ $\mathcal{N}\models \exists t<b\ \forall x\leq \mathrm{max}(%
\mathbf{v}_{m})\ \exists \left\langle r,i\right\rangle \in l_{s}\ \exists
z<t\ \mathrm{Sat}_{\Delta _{0}}(\delta _{r}(i,\mathbf{v}_{m},x,z)).$%
\smallskip

\noindent In particular, by choosing $x=d_{k}$ we obtain:\smallskip

\noindent $(xii)$ $\mathcal{N}\models \exists \left\langle r,i\right\rangle
<\,b\ \exists z<b(\left\langle r,i\right\rangle \in l_{s}\wedge \mathrm{Sat}%
_{\Delta _{0}}(\delta _{r}(i,\mathbf{v}_{m},d_{k},z)))$,\smallskip

\noindent which contradicts the definition of $l_{s}$.\smallskip

\noindent By $\Pi _{1}$-Overspill in $\mathcal{M}$ there is some $q\in
M\backslash I$ such that:\smallskip

\noindent $(xiii)$ $\mathcal{M}\models \exists x\leq \mathrm{max}(\mathbf{u}%
_{m})\ \forall \left\langle r,i\right\rangle \in l_{q}\ \forall z\ \lnot
\mathrm{Sat}_{\Delta _{0}}(\delta _{r},\mathbf{u}_{m},x,z).$\smallskip

\noindent Let $c$ be a witness in $\mathcal{M}$ to the $\exists x$ assertion
in $(xiii)$, and let $\mathbf{u}_{m+1}=\left( \mathbf{u}_{m},c\right) $ and $%
\mathbf{v}_{m+1}=\left( \mathbf{v}_{m},d_{k}\right) .$ It is easy to see
using $(xiii)$ that $(\ast _{m+1})$ holds with these choices of $\mathbf{u}%
_{m+1}$ and $\mathbf{v}_{m+1}.\medskip $

\noindent This concludes the recursive construction of $\left( u_{k}:k\in
\omega \right) $ and $\left( v_{k}:k\in \omega \right) $ satisfying
properties (I) through (IV).\hfill $\square $\textit{\medskip }

\noindent \textbf{3.3.1.}~\textbf{Corollary.}~\emph{Let} $\mathcal{M}$ \emph{%
and} $\mathcal{N}$ \emph{be countable nonstandard models of} $\mathrm{I}%
\Sigma _{1}$, \textit{and} $I$ \emph{be a proper cut shared by} $\mathcal{M}$
\emph{and} $\mathcal{N}$ \emph{that is closed under exponentiation. The
following are equivalent}:\textit{\smallskip }

\noindent \textbf{(1)}\textit{\ There is a proper initial embedding }$j$%
\textit{\ of }$\mathcal{M}$ \textit{into} $\mathcal{N}$ \textit{such that} $%
j(i)=i$ \emph{for all }$i\in I.$\textit{\smallskip }

\noindent \textbf{(2) }$\mathrm{Th}_{\Sigma _{1}}(\mathcal{M},i)_{i\in
I}\subseteq \mathrm{Th}_{\Sigma _{1}}(\mathcal{N},i)_{i\in I}$ \textit{and} $%
\mathrm{SSy}_{I}(\mathcal{M})=\mathrm{SSy}_{I}(\mathcal{N}).$\textit{%
\medskip }

\noindent \textbf{Proof.}$~(1)\Rightarrow (2)$ is again the easy direction.
To show that $(2)\Rightarrow (1)\mathbf{,}$ by Theorem 3.3 it suffices to
show (2) implies that there are $c\in M$ and $a,b\in N$ such that\emph{\ }%
for all\textrm{\ }$i\in I$ and\emph{\ }$\Delta _{0}$-formulae $\delta
(x,y,z) $,\emph{\ }if $\mathcal{M}\models \exists z\ \delta (i,c,z)$,\ then $%
\mathcal{N}\models \exists z<b\ \delta (i,a,z).$ Let $a=c=0$. We need to
show that for some $b\in N$ such that for all\textrm{\ }$i\in I$ and\emph{\ }%
$\Delta _{0}$-formulae $\delta (x,z)$,\emph{\ }if $\mathcal{M}\models
\exists z\ \delta (i,z)$,\emph{\ }then $\mathcal{N}\models \exists z<b\
\delta (i,z).$ Let\textbf{\ }$\left\langle \delta _{i}:i\in N\right\rangle $
be a canonical enumeration within $\mathcal{N}$ of all $\Delta _{0}$%
-formulae, and for $s\in N$ let $\varphi (s)$ be the following statement:

\begin{center}
$\exists y_{s}\ \forall \left\langle r,i\right\rangle <s\ \left[ \exists x\
\mathrm{Sat}_{\Delta _{0}}(\delta _{r}(i,x))\rightarrow \exists x<y_{s}\
\mathrm{Sat}_{\Delta _{0}}(\delta _{r}(i,x))\right] .$
\end{center}

\noindent By Strong $\Sigma _{1}$-collection in $\mathcal{N}$, $\varphi (s)$
holds in $\mathcal{N}$ for any $s\in N.$ In particular, if $s\in N\backslash
I$ then $y_{s}$ serves as our desired $b$. \hfill $\square $\textit{\medskip
}

\noindent \textbf{3.3.2.}~\textbf{Remark.}~For any element $a_{0}$ of $%
\mathcal{M}\models \mathrm{I}\Sigma _{1}$, let $\left( a_{n}:n<\omega
\right) $ be given by $\mathcal{M}\models a_{n+1}=2^{a_{n}}$; and consider:

\begin{center}
$I:=\{m\in M:$ $\exists n\in \omega $ such that $m<a_{n}\}.$
\end{center}

\noindent $I$ is by design closed under exponentiation; it also forms a
proper cut in $\mathcal{M}$ (thanks to the totality of the superexponential
function in $\mathcal{M}$). This makes it clear that the $(2)\Rightarrow (1)$
direction of Corollary 3.3.1 implies the $(1)\Rightarrow (2)$ direction of
Theorem 2.9.\textit{\medskip }

\noindent \textbf{3.4.~Theorem.}~\textit{For any countable nonstandard model}%
\textbf{\ }$\mathcal{M}$ \emph{of} $\mathrm{PA}$ \textit{there is a
composition preserving embedding }$j\longmapsto \widehat{j}$ \textit{of} $%
\mathrm{PISE}(\mathbb{Q})$ \textit{into} $\mathrm{PISE}(\mathcal{M})\mathrm{,%
}$ \textit{where} $\mathbb{Q}$\ \textit{is the ordered set of rationals.
Moreover, if} $j$ \textit{is contractive, then so is} $\widehat{j}.\medskip $

\noindent \textbf{Proof.}$~$Given a countable model $\mathcal{M}$ of $%
\mathrm{PA}$, choose $\mathcal{A}$ be the collection of subsets of $M$ that
are parametrically definable in $\mathcal{M}$, and let $\mathcal{N}_{\mathbb{%
Q}}$ be an elementary end extension of $\mathcal{M}$ as in Theorem 2.6.
Since $\mathcal{M}$ and $\mathcal{N}_{\mathbb{Q}}$ share the same standard
system and the same first order theory, Theorem 2.8 assures us that there is
a proper initial embedding $k:\mathcal{M}\rightarrow \mathcal{N}_{\mathbb{Q}%
} $ such that $M\subsetneq k(M).$ Let $M^{\ast }=k(M)$. By part (c) of
Theorem 2.6 we may choose $c_{q_{0}}\in $ $M^{\ast }\backslash M.$ Let $j\in
\mathrm{PISE}(\mathbb{Q})$ such that $j(\mathbb{Q})<q_{0}.$ By parts (d) and
(g) of Theorem 2.6:

\begin{center}
$\widehat{j}\in \mathrm{PISE}(\mathcal{N}_{\mathbb{Q}})$ and $\widehat{j}(N_{%
\mathbb{Q}})<c_{q_{0}}.$
\end{center}

\noindent Therefore $\widehat{j}(M^{\ast })<c_{q_{0}}\in M^{\ast }.$ Let $%
\widehat{j}_{M^{\ast }}$ be the restriction of $\widehat{j}$ to $M^{\ast }$.
Then $\widehat{j}_{M^{\ast }}\in \mathrm{PISE}(\mathcal{M}^{\ast })$ and the
desired embedding of $\mathrm{PISE}(\mathbb{Q})$ into $\mathrm{PISE}(%
\mathcal{M})$ is $j\mapsto k^{-1}\circ \widehat{j}_{M^{\ast }}\circ k$.
\hfill $\square $\textit{\medskip }

\noindent \textbf{3.4.1.}~\textbf{Remark.}$~$It is easy to see, using
Cantor's theorem asserting that any countable dense linear order without
endpoints is isomorphic to $\mathbb{Q}$, that $\mathbb{Q}$ carries a proper
initial self-embedding that is contractive. \textit{\medskip }

\noindent \textbf{3.4.2.}~\textbf{Corollary.}$~$\textit{Every countable
nonstandard model of }{\normalsize PA }\textit{carries a contractive proper
initial self-embedding}.\medskip

\noindent \textbf{Proof.}$~$Put Theorem 3.4 together with Remark
3.4.1.\hfill $\square $\textit{\medskip }

\noindent \textbf{3.4.3.}~\textbf{Proposition.}$~$\textit{For every countable%
} \textit{linear order }$\mathbb{L}$\textit{, there is a composition
preserving embedding }$j\mapsto \widehat{j}$ \textit{of} $\mathrm{SE}(%
\mathbb{L})$\textit{\ into} $\mathrm{SE}(\mathbb{Q})$. \textit{Moreover}:
\smallskip

\noindent \textbf{(a) }$\widehat{j}\in \mathrm{Aut}(\mathbb{Q})$ \textit{iff}
$j\in $ $\mathrm{Aut}(\mathbb{L)}$.\textit{\smallskip }

\noindent \textbf{(b) }$\widehat{j}\in \mathrm{ISE}(\mathbb{Q})$ \textit{iff}
$j\in $ $\mathrm{ISE}(\mathbb{L)}.$\textit{\smallskip }

\noindent \textbf{(c) }$\widehat{j}\in \mathrm{PISE}(\mathbb{Q})$ iff $j\in $
$\mathrm{PISE}(\mathbb{L)}$. \medskip

\noindent \textbf{Proof.}$~$Given a linear order $\mathbb{L}$, let $\mathbb{L%
}\times \mathbb{Q}$ be the lexicographic product of $\mathbb{L}$ and $%
\mathbb{Q}$ (intuitively $\mathbb{L}\times \mathbb{Q}$ is the result of
replacing each point in $\mathbb{L}$ by a copy of $\mathbb{Q}$). $\mathbb{L}%
\times \mathbb{Q}$ is clearly a countable dense linear order with no end
points. Therefore when $\mathbb{L}$\ is countable, $\mathbb{L}\times \mathbb{%
Q}$ is isomorphic to $\mathbb{Q}$ by Cantor's theorem mentioned in Remark
3.4.1. So it suffices to find a composition preserving embedding of $\mathrm{%
SE}(\mathbb{L})$\textit{\ }into $\mathrm{SE}(\mathbb{L}\times \mathbb{Q})$
that satisfies (a), (b), and (c). Given $j\in \mathrm{SE}(\mathbb{L}),$ let $%
\widehat{j}:\mathbb{L}\times \mathbb{Q\rightarrow L}\times \mathbb{Q}$ by $%
\widehat{j}(l,q)=(j(l),q).$ A routine reasoning shows that $\widehat{j}\in
\mathrm{SE}(\mathbb{L}\times \mathbb{Q}),$ and the embedding $j\mapsto
\widehat{j}$ is composition preserving. Properties (a), (b), and (c) are
equally easy to verify.\hfill $\square \medskip $

\noindent \textbf{3.4.4.}~\textbf{Remark.}$~$Let $\mathcal{M}=(M,<,\cdot
\cdot \cdot )$ be a linearly ordered structure. $\mathrm{SE}(\mathcal{M})$
is a sub-semigroup of \textrm{SE}$(M,<),$ therefore by Proposition 3.4.3 $%
\mathrm{SE}(\mathcal{M})$\ is embeddable into $\mathrm{SE}(\mathbb{Q})$; $%
\mathrm{Aut}(\mathcal{M})$\ is embeddable in $\mathrm{Aut}(\mathbb{Q})$; $%
\mathrm{ISE}(\mathcal{M})$\ is embeddable in $\mathrm{ISE}(\mathbb{Q})$; and
$\mathrm{PISE}(\mathcal{M})$ is embeddable in $\mathrm{PISE}(\mathbb{Q}).$%
\texttt{\bigskip }

\begin{center}
\textbf{4.~THE LONGEST INITIAL SEGMENT OF FIXED POINTS}\texttt{\bigskip }
\end{center}

In this section we establish the first principal result of this paper
(Theorem 4.1) by an elaboration of the back-and-forth proof of Theorem 3.3.
The $(2)\Rightarrow (3)$ direction of Theorem 4.1 fine-tunes the $%
(1)\Rightarrow (2)$ direction of Theorem 2.9, since as pointed out in Remark
3.3.2 proper cuts closed under exponentiation can be found arbitrarily high
in every nonstandard model of $\mathrm{I}\Sigma _{1}$. \medskip

\noindent \textbf{4.1.}$~$\textbf{Theorem}.$~$\textit{Suppose }$I$\textit{\
is a proper initial segment of a countable nonstandard model }$\mathcal{M}$
\textit{of} $\mathrm{I}\Sigma _{1}$\textit{. The following are equivalent}%
:\smallskip

\noindent \textbf{(1)} $I=$ $\mathrm{I}_{\mathrm{fix}}(j)$ \textit{for some}
\textit{self-embedding }$j$ \textit{of} $\mathcal{M}$\textit{.}\smallskip

\noindent \textbf{(2)} $I$\textit{\ is closed under exponentiation.}%
\smallskip

\noindent \textbf{(3)} $I=$ $\mathrm{I}_{\mathrm{fix}}(j)$ \textit{for some}
\textit{proper} \textit{initial self-embedding }$j$ \textit{of} $\mathcal{M}$%
\textit{.}\medskip

\noindent \textbf{Proof.}$~(1)\Rightarrow (2)$ follows immediately from
Lemma 3.1.1 and Theorem 3.2; and $(3)\Rightarrow (1)$ is trivial; so it
suffices to establish $(2)\Rightarrow (3)$. By the proof of Corollary 3.3.1
we can let $a=c=0$, and let $b$ be a large enough element of $\mathcal{M}$
such that for all\ $i\in I$ and\emph{\ }all\emph{\ }$\Delta _{0}$-formulae $%
\delta (x,y,z)$ we have:

\begin{center}
$\mathcal{M}\models \exists z\ \delta (i,c,z)\rightarrow \exists z<b\ \delta
(i,a,z)$.
\end{center}

\noindent Assume (2). In order to produce the desired embedding $j$
satisfying (3) we will elaborate the proof of Theorem 3.3 by adding a third
layer of recursion to the proof of Theorem 3.3. More specifically, at stage $%
m=3k$ we will do the same as stage $m=2k$ of the proof of Theorem 3.3, and
at stage $m=3k+1$ we will do the same as stage $m=2k+1$ of the proof of
Theorem 3.3. In order to describe the construction for stages $m$ of the
form $3k+2$, we first establish the following lemma:\medskip

\noindent \textbf{4.1.1.}$~$\textbf{Lemma}$~$\textit{Suppose} $\mathbf{u}$
\textit{and} $\mathbf{v}$ \textit{are finite tuples of the same length from}
$\mathcal{M}$ \textit{that satisfy}:\smallskip

\noindent \textbf{(I)}$\ \mathcal{M}\models \exists z\ \delta (i,\mathbf{u}%
,z)\rightarrow \exists z<b\ \delta (i,\mathbf{v},z)$ \textit{for any} $i\in
I $ \textit{and any }$\delta (x,\mathbf{y},z)\in \Delta _{0}.\smallskip $

\noindent \textit{Then} \textit{for any} $d\in M\backslash I$ \textit{there
are \textbf{distinct} }$u,v\in M$ \textit{such that} $u<d$ \textit{and}%
:\smallskip

\noindent \textbf{(II)} $\mathcal{M}\models \exists z\ \delta (i,\mathbf{u}%
,u,z)\rightarrow \exists z<b\ \delta (i,\mathbf{v},v,z)$ \textit{for any} $%
i\in I$ \textit{and any }$\delta (x,\mathbf{y},w,z)\in \Delta _{0}$.\medskip

\noindent \textbf{Proof}.$~$Assume (I) holds and suppose $d\in M\backslash
I. $ Let\textbf{\ }$\left\langle \delta _{i}:i\in M\right\rangle $ be a
canonical enumeration within $\mathcal{M}$ of all $\Delta _{0}$-formulae.
For $s\in I$ and $x<d$, let:

\begin{center}
$H_{s,x}:=\left\{ \left\langle r,i\right\rangle <s:\exists z\ \mathrm{Sat}%
_{\Delta _{0}}(\delta _{r}(i,\mathbf{u},x,z))\right\} $.
\end{center}

\noindent Then define $f_{s}:\underline{d}\rightarrow \underline{2^{s+1}}$
in $\mathcal{M}$ for $x<d$ via:

\begin{center}
$f_{s}(x)=\sum\limits_{\left\langle r,i\right\rangle \in
H_{s,x}}2^{\left\langle r,i\right\rangle }$.
\end{center}

\noindent Note that $f_{s}(x)\leq \sum\limits_{k<s}2^{k}=2^{s+1}-1$, which
coupled with the closure of $I$ under exponentiation implies:

\begin{center}
$d>2^{s+1}>$ $f_{s}(x)$.
\end{center}

\noindent On the other hand, for each $x$, $f_{s}(x)$ is $\Sigma _{1}$%
-minimal (in parameters $x$ and $s$), and therefore the graph of $f_{s}$ is $%
\Delta _{0}(\Sigma _{1})$-definable in $\mathcal{M}$, so by $\Delta
_{0}(\Sigma _{1})$-PHP, $f_{s}$ is not one-to-one, and we may therefore
choose distinct $u,u^{\prime }<d$ such that $f_{s}(u)=f_{s}(u^{\prime })$.
Let $\varphi (s)$ be the formula:

\begin{center}
$\exists u,u^{\prime }<d\ \left( (u\neq u^{\prime })\wedge \theta
(s,u,u^{\prime })\right) ,$
\end{center}

\noindent where $\theta (s,u,u^{\prime })$ is:

\begin{center}
$\forall \left\langle r,i\right\rangle <s\ \left[ \exists z\ \mathrm{Sat}%
_{\Delta _{0}}\left( \delta _{r}(i,\mathbf{u},u,z)\right) \leftrightarrow
\exists z\ \mathrm{Sat}_{\Delta _{0}}\left( \delta _{r}(i,\mathbf{u}%
,u^{\prime },z)\right) \right] .$
\end{center}

\noindent The definition of $f_{s}$ makes it evident that $\mathcal{M}%
\models \varphi (s)$ for each $s\in I$. Since $\varphi (s)$ is a $\Delta
_{0}(\Sigma _{1})$ statement, by $\Delta _{0}(\Sigma _{1})$-Overspill in $%
\mathcal{M}$ there is some $p\in M\backslash I$ such that $\mathcal{M}%
\models \varphi (p).$ Therefore there are distinct $u,u^{\prime }<d$ such
that for each $i\in I$ and each $\Delta _{0}$-formula $\delta $ we
have:\smallskip

\noindent $(i)$ $\mathcal{M}\models \exists z\ \delta (i,\mathbf{u}%
,u,z)\leftrightarrow \exists z\ \delta (i,\mathbf{u},u^{\prime },z)$%
.\smallskip

\noindent On the other hand, by the proof of the `forth' direction (the $%
m=2k $ case) of Theorem 3.3, we can find distinct $w$ and $w^{\prime }$ such
that the following holds for each $\Delta _{0}$-formula $\delta $:\smallskip

\noindent $(ii)$ $\mathcal{M}\models \exists z\ \delta (i,\mathbf{u}%
,u,u^{\prime },z)\rightarrow \exists z<b\ \delta (i,\mathbf{v},w,w^{\prime
},z)$.\smallskip

\noindent Since at least one of the two statements $\{u\neq w$, $u\neq
w^{\prime }\}$ is true, we can choose $v\in \{w,w^{\prime }\}$ such that $%
u\neq v.$ It is easy to see using $(i)$ and $(ii)$ that this choice of $u$
and $v$ satisfy (II). \hfill $\square $ Lemma 4.1.1\medskip

Fix a sequence $\left( d_{k}:k\in \omega \right) $ that is downward cofinal
in $M\setminus I.$ Suppose $m=3k+2$ and we have $\mathbf{u}_{m}$ and $%
\mathbf{v}_{m}$ satisfying condition $(\ast _{m})$ of the proof of Theorem
3.3 for $\mathcal{N}:=\mathcal{M}$. Apply Lemma 4.1.1 with $\mathbf{u}:=%
\mathbf{u}_{m}$, $\mathbf{v}:=\mathbf{v}_{m},$ and $d:=d_{k}$ to get hold of
$u$ and $v$ satisfying (II) of Lemma 4.1.1; and then we define $\mathbf{u}%
_{m+1}:=\left( \mathbf{u}_{m},u\right) $ and $\mathbf{v}_{m+1}:=\left(
\mathbf{v}_{m},v\right) .$ This makes it clear that the proper initial
self-embedding $j$ of $\mathcal{M}$ defined by $j(u_{k})=v_{k}$ fixes each $%
i\in I$ but moves elements arbitrarily low in $M\backslash I$.\hfill $%
\square \medskip $

\noindent \textbf{4.1.2.}$~$\textbf{Remark}.$~$By Remark 3.3.2 there are
unboundedly many cuts in a nonstandard model of I$\Sigma _{1}$ that are
closed under exponentiation. Therefore Theorem 3.3 is a strengthening of the
$(1)\Rightarrow (2)$ direction of Theorem 2.9. Also, it is easy to see
(using an overspill argument) that in nonstandard models of $\mathrm{I}%
\Delta _{0}$ cuts that are closed under exponentiation can be found
arbitrarily low in the nonstandard part of $\mathcal{M}$.\textbf{\bigskip }

\begin{center}
\textbf{5.~FIXED POINT SETS THAT ARE INITIAL SEGMENTS\bigskip }
\end{center}

This section is devoted to the second main result of this paper (Theorem
5.1). See also Remark 5.1.1.\medskip

\noindent \textbf{5.1.}$~$\textbf{Theorem}.$~$\textit{Suppose }$I$\textit{\
is a proper initial segment of a countable nonstandard model }$\mathcal{M}$%
\textit{\ of }$\mathrm{I}\Sigma _{1}$\textit{. The following are equivalent}%
:\smallskip

\noindent \textbf{(1)} $I=\mathrm{Fix}(j)$ \textit{for some self-embedding }$%
j$ \textit{of} $\mathcal{M}.$\smallskip

\noindent \textbf{(2)} $I$\textit{\ is a strong cut of }$\mathcal{M}$\textit{%
, and }$I\prec _{\Sigma _{1}}\mathcal{M}.$\smallskip

\noindent \textbf{(3)} $I=\mathrm{Fix}(j)$ \textit{for some proper initial
self-embedding }$j$ \textit{of} $\mathcal{M}.$\textit{\medskip }

\noindent \textbf{Proof.}$~$Since\textbf{\ }$(3)\Rightarrow (1)$ is trivial,
it suffices to show $(1)\Rightarrow (2)$ and $(2)\Rightarrow (3).\medskip $

\noindent To verify $(1)\Rightarrow (2)\mathbf{,}$ suppose (1) holds and let
$\widetilde{f}\in M$ code an $\mathcal{M}$-finite function $f$ whose domain
includes $I$. It is easy to see that $\widetilde{f}\notin I.$ So if $%
\widetilde{g}:=j(\widetilde{f})$, then $\widetilde{g}\notin I,$ and $%
\widetilde{f}\neq \widetilde{g}$. Therefore, in light of the assumption that
$I=\mathrm{Fix}(j),$ if $g$ is the function that is coded by $\widetilde{g}$%
, then:

\begin{center}
$\forall i\in I$ $[f(i)=g(i)\Longleftrightarrow f(i)\in I].$
\end{center}

\noindent We wish to find $s\in M\backslash I$ such that for all $i\in I,$ $%
f(i)\notin I$\ iff $s<f(i).$ Fix $d\in M$ such that $I<d$ and the interval $%
[0,d]\subseteq \mathrm{dom}(f)\cap \mathrm{dom}(g).$ Without loss of
generality there is some $i_{0}\in I$ with $f(i_{0})\notin I.$ Consider the
function $h(x)$ defined within $\mathcal{M}$ on the interval $[i_{0},d]$ by:

\begin{center}
$h(x):=\mu y\leq d$\ $[\exists z\leq x(y=f(z)\neq g(z)],$
\end{center}

\noindent where $\mu y\leq d$ is the modified least search operator, defined
via the following:

\begin{center}
$[z:=\mu y\leq d$ $\varphi (y)]$

iff

[$z$ is the first $y$ such that $\varphi (y)$, if $\exists y\leq d\ \varphi
(y);$ else $z=d$].\smallskip
\end{center}

\noindent Note that if $i\in I$, then $h(i)\notin I$, and if $i_{0}\leq
i\leq i^{\prime }$, then $h(i^{\prime })\leq h(i)$. Moreover:\smallskip

\noindent $(i)$ The graph of $h$ is defined by a $\Delta _{0}$-formula $%
\varphi (x,y)$ with parameters $\widetilde{f}$, $\widetilde{g}$, and $d$%
.\smallskip

\noindent $(ii)$ $i<h(i)$ for all $i\in I$ such that $i\geq i_{0}.$\smallskip

\noindent Therefore, by putting $(i)$ together with $(ii)$ and $\Delta _{0}$%
-Overspill we may conclude that there is some $s\in M\backslash I$ such that
$s<h(s)$ holds in $\mathcal{M}.$ This shows that $s$ is the desired lower
bound for elements of the form $f(i)$, where $i\in I$ and $f(i)\notin I$.
This concludes the verification that $I$\ is a strong cut of $\mathcal{M}$.
On the other hand, since we are assuming that (1) holds, Theorem 3.1 assures
that $I\prec _{\Sigma _{1}}\mathcal{M}$, so (2) holds.\medskip

\noindent To establish $(2)\Rightarrow (3),$ suppose (2) holds. We first
note that by Theorem 2.5, $(I,\mathrm{SSy}_{I}(\mathcal{M}))\models \mathrm{%
ACA}_{0}$. By Theorem 2.6 we can build $\mathcal{N}_{\mathbb{Q}}\succ _{%
\mathrm{end}}I$ (where $\mathbb{Q}$\ is the ordered set of rationals) such
that:\smallskip

\noindent $(iii)$ $\mathrm{SSy}_{I}(\mathcal{M)}=\mathrm{SSy}_{I}(\mathcal{N}%
_{\mathbb{Q}}\mathcal{)}$, and \smallskip

\noindent $(iv)$ $\mathbb{Q}^{\prime }:=\{c_{q}:q\in \mathbb{Q}\}$ is an
isomorphic copy of $\mathbb{Q}$ and is downward cofinal in $\mathcal{N}_{%
\mathbb{Q}}\backslash I$. \smallskip

\noindent On the other hand, since $I\prec _{\Sigma _{1}}\mathcal{M}$ we may
infer that $\mathrm{Th}_{\Sigma _{1}}(\mathcal{M},i)_{i\in I}=\mathrm{Th}%
_{\Sigma _{1}}(\mathcal{N}_{\mathbb{Q}},i)_{i\in I}$, which together with $%
(iii)$\ and Corollary 3.3.1 allows us to get hold of an initial embedding $k:%
\mathcal{M}\rightarrow \mathcal{N}_{\mathbb{Q}}$ such that $k$ pointwise
fixes each $i\in I.$ Let $M^{\ast }$ be the range of $k$. By $(iv)$ there is
some $q_{0}\in \mathbb{Q}$ and $m_{0}^{\ast }\in M^{\ast }$ such
that:\smallskip

\noindent $(v)\ c_{q_{0}}<m_{0}^{\ast }.$ \smallskip

\noindent Let $j_{0}:\mathbb{Q}\rightarrow \mathbb{Q}$ be a proper initial
self-embedding of $\mathbb{Q}$ whose range is bounded above by $q_{0}.$ By
Theorem 2.6 the range of the induced initial self-embedding $\widehat{j}_{0}$
of $\mathcal{N}_{\mathbb{Q}}$ is bounded above by $c_{q_{0}}$ and $\mathrm{%
Fix}(\widehat{j}_{0})=I$. Coupled with $(v)$ this shows that $\widehat{j_{0}}%
(M^{\ast })\subsetneq M^{\ast }.$ So we can identify $\mathcal{M}$ with its
isomorphic copy $\mathcal{M}^{\ast }$ to complete the proof; in other words
the desired $j\in \mathrm{PISE}(\mathcal{M})$ such that $\mathrm{Fix}(j)=I$
is given by $j:=k^{-1}\widehat{j}_{0}k.$\hfill $\square \medskip $

\noindent \textbf{5.1.1.}$~$\textbf{Remark}.$~$For each $n\in \omega $,
there is a countable model of $\mathrm{I}\Sigma _{n}$ which does not carry a
proper cut $I$ satisfying (2) of Theorem 5.1. To see this, first note that
(2) implies that $\mathcal{M}$ $\models \mathrm{Con}(\mathrm{I}\Sigma _{n})$
for each $n<\omega $ since $\mathrm{PA}$ holds in $I$ by Theorem 2.5, $%
\mathrm{Con}(\mathrm{I}\Sigma _{n})$ is a $\Pi _{1}$-statement, and it is
well-known \cite[Ex.~10.8]{Kaye's text} that $\mathrm{Con}(\mathrm{I}\Sigma
_{n})$ is provable in $\mathrm{I}\Sigma _{n+1}$ for each $n\in \omega $. On
the other hand, $\mathrm{Con}(\mathrm{I}\Sigma _{n})$ is unprovable in $%
\mathrm{I}\Sigma _{n}$ by G\"{o}del's second incompleteness theorem, and
therefore there is a countable nonstandard model\textit{\ }$\mathcal{M}_{0}$
of\textit{\ }$\mathrm{I}\Sigma _{n}+\lnot \mathrm{Con}(\mathrm{I}\Sigma
_{n}).$ Such a model $\mathcal{M}_{0}$\textit{\ }has no cut\textit{\ }that
satisfies condition $(2)$\ of Theorem 5.1. However, if $\mathcal{M}$ is a
countable nonstandard model of $\mathrm{PA}$, then by using a variation of
the proof of Tanaka's theorem in \cite{Me Tanaka}, for any $n\in \omega $ we
can find a strong cut $I$\ arbitrarily high in $\mathcal{M}$ such that $%
I\prec _{\Sigma _{n}}\mathcal{M}.$ Tin Lok Wong has also pointed out to us
that there are countable models $\mathcal{M}_{0}$ of $\mathrm{I}\Sigma _{1}$
in which there is no proper cut $I$ such that $I\prec _{\Sigma _{1}}\mathcal{%
M}_{0}$. Such a model $\mathcal{M}_{0}$ can be readily obtained by choosing $%
\mathcal{M}_{0}$ as $H^{1}(\mathcal{M}),$ where $\mathcal{M}\models \mathrm{I%
}\Sigma _{1}$ and $H^{1}(\mathcal{M})$ is defined as in \cite[Ch.~IV,
Def.~1.32]{Hajek-Pudlak text}. \textbf{\bigskip }

\begin{center}
\textbf{6.~MINIMAL FIXED POINTS\bigskip }
\end{center}

In this section we establish our final principal result (Theorem 6.1). The
proof of Theorem 6.1 is rather complex and based on several technical
lemmas, which were inspired by, and can be seen as miniaturized analogues of
Lemmas 8.6.4, 8.6.6, and 8.6.2 of \cite{Kossak-Schmerl} (which were
originally established in the joint work of Kaye, Kossak, and Kotlarski \cite%
{Richard-Roman-Henryk}). $\medskip $

Recall from Theorem 3.1 that $K^{1}(\mathcal{M})\subseteq \mathrm{Fix}(j)$
for every $j\in \mathrm{SE}(\mathcal{M}),$ where $\mathcal{M}\models \mathrm{%
I}\Delta _{0}+\mathrm{Exp}$. It is also straightforward to modify the proof
of the basic Friedman embedding theorem \cite[Thm.~12.3]{Kaye's text} to
show that if $\mathcal{M}$ is a countable nonstandard model of $\mathrm{%
I\Sigma }_{1},$ and $m\in M\backslash K^{1}(\mathcal{M})$, then there is
some $j\in \mathrm{PISE}(\mathcal{M})$ such that $j(m)\neq m.$ These results
motivate the question whether every countable nonstandard model $\mathcal{M}%
\models \mathrm{I}\Sigma _{1}$ has a proper initial self-embedding that
moves \textit{all} elements of $M\backslash K^{1}(\mathcal{M}).$ Theorem 6.1
provides a complete answer to this question.\medskip

\noindent \textbf{6.1.}$~$\textbf{Theorem}.$~$\textit{The following are
equivalent for a countable nonstandard model }$\mathcal{M}$\textit{\ of }$%
\mathrm{I}\Sigma _{1}$:\smallskip

\noindent \textbf{(1)} $\mathrm{Fix}(j)=K^{1}(\mathcal{M})$ \textit{for some
self-embedding }$j$ \textit{of} $\mathcal{M}.$\smallskip

\noindent \textbf{(2)} $\mathbb{N}$ \textit{is a strong cut of} $\mathcal{M}$%
.\smallskip

\noindent \textbf{(3)} $\mathrm{Fix}(j)=K^{1}(\mathcal{M})$ \textit{for some
proper initial self-embedding }$j$ \textit{of} $\mathcal{M}.$\textit{%
\medskip }

\noindent \textbf{Proof.}$~$Since $(3)\Rightarrow (1)$ is trivial, it
suffices to show that $(1)\Rightarrow (2)$, and $(2)\Rightarrow (3).\medskip
$

\begin{center}
\textbf{Proof} \textbf{of} $\mathbf{(1)\Rightarrow (2)}$ \textbf{of Theorem
6.1} $\medskip $
\end{center}

The proof is based on Lemma 6.1.1, 6.1.2, and 6.1.4 below.\textit{\medskip }

\noindent \textbf{6.1.1.}$~$\textbf{Lemma.}$~$\textit{If }$\mathbb{N}$%
\textit{\ is not a strong cut of }$\mathcal{M}\models \mathrm{I}\Delta _{0}$%
\textit{, then for any self-embedding }$j$\textit{\ of }$\mathcal{M}$\textit{%
, the nonstandard fixed points of }$j$\textit{\ are downward cofinal in the
nonstandard part of }$\mathcal{M}$. \textit{\medskip }

\noindent \textbf{Proof.}~Suppose that $\mathbb{N}$ is not strong in $%
\mathcal{M}$. Then there is some function $f$ coded in $\mathcal{M}$ whose
domain is of the form $\underline{c}$ for some nonstandard $c$, and such
that $D:=\left\{ f(n):n\in \mathbb{N\ }\mathrm{and}\ f(n)\in M\backslash
\mathbb{N}\right\} $ is downward cofinal in the nonstandard part of $%
\mathcal{M}$. Let $j$ be a self-embedding of $\mathcal{M}$, and let $g:=j(f)$%
. We observe that for each standard number $n$ the statement $P(n)$ holds in
$\mathcal{M}$, where:

\begin{center}
$P(z):=$ \textquotedblleft For all $x,y<z,\ f(x)=y$ iff $g(x)=y$%
\textquotedblright .
\end{center}

\noindent Since $P(z)$ is a $\Delta _{0}$-formula (with parameters $f$ and $%
g $), by $\Delta _{0}$-Overspill for any nonstandard $k\in M$ there is some
nonstandard $c<k$ such that $P(c)$ holds in $\mathcal{M}$. So it suffices to
show that there is a nonstandard fixed point below any such $c$. Going back
to the set $D$ defined earlier, let $n_{0}\in \mathbb{N}$ such that $%
f(n_{0}) $ is nonstandard and $f(n_{0})<c$. Note that $f(n_{0})=g(n_{0})$
since $P(c)$ holds in $\mathcal{M}$, therefore:

\begin{center}
$j(f(n_{0}))=j(f)(j(n_{0}))=g(n_{0})=f(n_{0}).$
\end{center}

\noindent So $f(n_{0})$ is the desired nonstandard fixed point of $j$ below $%
c$.\hfill $\square $\textit{\medskip }

\noindent \textbf{6.1.2.}$~$\textbf{Lemma.}$~$\textit{Suppose }$\mathbb{N}$
\textit{is not a strong cut of} $\mathcal{M}\models \mathrm{I}\Sigma _{1}$.
\textit{Then for every element }$a\in M$\textit{\ and any self-embedding }$j$%
\textit{\ of }$\mathcal{M}$ \textit{there is an element }$b\in \mathrm{Fix}%
(j)$\textit{\ such that}: \medskip

\begin{center}
$\mathrm{Th}_{\Sigma _{1}}(\mathcal{M},a)\subseteq \mathrm{Th}_{\Sigma _{1}}(%
\mathcal{M},b).$
\end{center}

\noindent \textbf{Proof.}$~$Let\textbf{\ }$\left\langle \sigma _{i}(x):i\in
M\right\rangle $ be a canonical enumeration within $\mathcal{M}$ of all $%
\Sigma _{1}$-formulae in one free variable $x,$ with $\sigma _{i}(x)=\exists
y\ \delta _{i}(x,y),$ where $\delta _{i}$ is a $\Delta _{0}$-formula in the
sense of $\mathcal{M}$. Recall that $(x)_{i}$ refers to the $i$-th
coordinate of the sequence canonically coded by $x$, and the graph of $%
(x)_{i}$ is $\Delta _{0}$-definable.\textit{\medskip }

Given $a\in M$, for any $k\in M,$ $\left\{ i<k:\exists y\ \mathrm{Sat}%
_{\Delta _{0}}\left( \delta _{i}(a,y)\right) \right\} $ is coded by some $%
\mathcal{M}$-finite $s_{k}$ thanks to part (b) of Theorem 2.3 and the fact
that $\mathrm{Sat}_{\Delta _{0}}\in \Sigma _{1}^{\mathrm{I\Sigma }_{1}}.\ $%
Note that the mapping $k\mapsto s_{k}$ is $\Sigma _{1}$-definable in $(%
\mathcal{M},a)$. This makes it clear that for any $k\in M$ there is $%
c_{k}\in M$ such that:

\begin{center}
$\mathcal{M}\models c_{k}=\min \left\{ m\in M:\mathrm{Sat}_{\Delta
_{0}}\left( \bigwedge\limits_{iEs_{k}}\delta _{i}(\left( m\right)
_{0},\left( m\right) _{i+1})\right) \right\} $
\end{center}

\noindent We observe that:\medskip

\noindent $(i)$ For each nonstandard $k\in M\quad \mathrm{Th}_{\Sigma _{1}}(%
\mathcal{M},a)\subseteq \mathrm{Th}_{\Sigma _{1}}\left( \mathcal{M},\left(
c_{k}\right) _{0}\right) .$\medskip

\noindent Fix a nonstandard $i\in M$ choose $d\in M$ with $\left( d\right)
_{k}=\left( c_{k}\right) _{0}$ for all $k<i$. Note that $\left( d\right)
_{n}\in K^{1}(\mathcal{M})$ for $n\in \omega $, and therefore $j(\left(
d\right) _{n})=\left( d\right) _{n}$ for $n\in \omega .$ On the other hand,
if we let $e:=j(d),$ then for $n\in \omega $: \medskip

\begin{center}
$j(\left( d\right) _{n})=\left( j\left( d\right) \right) _{j(n)}=\left(
e\right) _{n}.$
\end{center}

\noindent This shows that $\left( d\right) _{n}=\left( e\right) _{n}$ for $%
n\in \omega $, so if we let:

\begin{center}
$\varphi (x):=\forall i<x$ $\left( d\right) _{i}=\left( e\right) _{i},$
\end{center}

\noindent then $\varphi (n)$ holds in $\mathcal{M}$ for each $n\in \omega $;
hence by $\Delta _{0}$-Overspill there is some nonstandard $n^{\ast }$ below
$i$ such that $\left( d\right) _{k}=\left( e\right) _{k}$ for all $k\leq
n^{\ast }.$ Therefore by Lemma 6.1.1 there is a nonstandard $k\in M$ that is
below $n^{\ast }$ such that:\medskip

\noindent $(ii)$ $(d)_{k}=(e)_{k}$ and $j(k)=k.$ \medskip

\noindent Since $\left( d\right) _{k}=\left( c_{k}\right) _{0}$ by design,
in light of $(i)$ the proof of our lemma will be complete once we observe
that $(d)_{k}\in \mathrm{Fix}(j)$ since by $(ii)$ we have:

\begin{center}
$j((d)_{k})=\left( j\left( d\right) \right) _{j(k)}=(e)_{k}=(d)_{k}.$

\hfill $\square $\medskip
\end{center}

It is convenient to employ the notion of a \textit{partial recursive
function of} $\mathcal{M}$ in order to state the next lemma; this notion
will also play a key role in the proof of $(2)\Rightarrow (3)$ of Theorem
6.1.\medskip

\noindent \textbf{6.1.3.}$~$\textbf{Definition.}$~$A partial function $f$
from $M$ to $M$ is a \textit{partial recursive function} of $\mathcal{M}$
iff the graph of $f$ is definable in $\mathcal{M}$ by a parameter-free $%
\Sigma _{1}$-formula; i.e., there is some $\Delta _{0}$-formula $\delta
(x,y,z)$ such that for all elements $r$ and $s$ of $\mathcal{M}$:

\begin{center}
$f(r)=s$ iff $\mathcal{M}\models \exists z\ \delta (r,s,z)$.
\end{center}

\noindent Given such an $f$, we will write $\left[ f(x)\downarrow \right] $
as an abbreviation for $\exists y\exists z\ \delta (x,y,z),$ and $\left[
f(x)\downarrow \right] ^{<w}$ as an abbreviation for:

\begin{center}
$\exists y,z<w\ \delta (x,y,z).$
\end{center}

\noindent Note that a partial recursive function $f$ naturally induces for
each positive $n\in \omega $ a partial function from $M^{n}$ to $M$, which
we will also denote by $f,$ via:

\begin{center}
$f(a_{1},\cdot \cdot \cdot ,a_{n}):=f(\left\langle a_{1},\cdot \cdot \cdot
,a_{n}\right\rangle ).$
\end{center}

\begin{itemize}
\item We shall use $\mathcal{F}$ to denote the collection of all partial
recursive functions of $\mathcal{M}$.
\end{itemize}

\noindent \textbf{6.1.4.}$~$\textbf{Lemma. }\textit{If} $\mathcal{M}\models
\mathrm{I}\Delta _{0},$ \textit{then}:

\begin{center}
$K^{1}(\mathcal{M})=\left\{ f(0):f\in \mathcal{F\ }\mathrm{and\ }\mathcal{M}%
\models \left[ f(0)\downarrow \right] \right\} .$
\end{center}

\noindent \textbf{Proof.}$~$This is an immediate consequence of part (b) of
Lemma 3.1.2.\hfill $\square $\medskip

With the above lemmas in place we are now ready to present the proof of $%
(1)\Rightarrow (2)$ by demonstrating its contrapositive. Suppose $\mathbb{N}$
is not a strong cut of $\mathcal{M}$. Consider the type $p(x)$ consisting of
the $\Sigma _{1}$-formulae of the form $\left[ f(0)\downarrow \right] \wedge
x\neq f(0)$, as $f$ ranges over the partial recursive functions of $\mathcal{%
M}$. By Lemma 6.1.4 no element of $K^{1}(\mathcal{M})$ realizes $p(x)$, and
yet $p(x)$ is realized by every element of $M\backslash K^{1}(\mathcal{M})$,
and of course $M\backslash K^{1}(\mathcal{M})\neq \varnothing $ ($\mathrm{B}%
\Sigma _{1}$ holds in $\mathcal{M}$, but not in $K^{1}(\mathcal{M})$ by $n=0$
case of part (b) of Theorem 2.4). In particular, if $a$ is chosen as an
element of $M\backslash K^{1}(\mathcal{M})$ then for $b\in M,$ $\mathrm{Th}%
_{\Sigma _{1}}(\mathcal{M},a)\subseteq \mathrm{Th}_{\Sigma _{1}}\left(
\mathcal{M},b\right) $ implies $b\notin K^{1}(\mathcal{M})$. Hence $K^{1}(%
\mathcal{M})\neq \mathrm{Fix}(j)$ by Lemma 6.1.2. This concludes the proof
of $(1)\Rightarrow (2)$ of Theorem 6.1.\hfill $\square $\medskip

\begin{center}
\textbf{Proof of }$\mathbf{(2)\Rightarrow (3)}$ \textbf{of Theorem 6.1}%
\medskip
\end{center}

Assume (2). Since $\mathrm{Sat}_{\Sigma _{1}}$ has a $\Sigma _{1}$%
-description in $\mathcal{M}$ and strong $\Sigma _{1}$-collection holds in $%
\mathcal{M}$, there is a sufficiently large $b\in M$ such that: \smallskip

\noindent $(\triangledown )$ For all $\Delta _{0}$-formulae $\delta (x),$ $%
\mathcal{M}\models \exists x\ \delta (x)\rightarrow \exists x<b\ \delta (x).$%
\smallskip

\noindent\ Note that $(\triangledown )$ is equivalent to:\smallskip

\noindent $(\blacktriangledown )$ For all $f\in \mathcal{F},$ $\mathcal{M}%
\models \left[ f(0)\downarrow \right] \rightarrow \left[ f(0)\downarrow %
\right] ^{<b}.$\smallskip

\noindent It is clear that the proof of (3) will be complete by setting $%
j(u_{k})=v_{k}$ once we have two sequences $\left( u_{r}:r<\omega \right) $
and $\left( v_{r}:r<\omega \right) $ that satisfy the following four
conditions:\medskip

\noindent (I) $M=\left\{ u_{r}:r<\omega \right\} .\smallskip $

\noindent (II) $\left\{ v_{r}:r<\omega \right\} $ is an initial segment of $%
\mathcal{N}$, and each $v_{r}<b.\smallskip $

\noindent (III) For each positive $n<\omega $, the following two properties $%
P(\mathbf{u},\mathbf{v})$ and $Q(\mathbf{u},\mathbf{v})$ hold for $\mathbf{u}%
=\left\langle u_{r}:r<n\right\rangle $, and $\ \mathbf{v}=\left\langle
v_{r}:r<n\right\rangle $:

\begin{center}
$P(\mathbf{u},\mathbf{v})$: For every $f\in \mathcal{F}$, $\mathcal{M}%
\models \left[ f(\mathbf{u})\downarrow \right] \rightarrow \left[ f(\mathbf{v%
})\downarrow \right] ^{<b}$. \medskip

$Q(\mathbf{u},\mathbf{v})$: For every $f\in \mathcal{F}$, if $\mathcal{M}%
\models \left[ f(\mathbf{u})\downarrow \right] $ and $f(\mathbf{u})\notin
K^{1}(\mathcal{M}),$ then $\mathcal{M}\models \left[ f(\mathbf{v})\downarrow %
\right] ^{<b}\wedge f(\mathbf{u})\neq f(\mathbf{v}).$
\end{center}

\noindent Note that $P(\mathbf{u},\mathbf{v})$ is equivalent to asserting
that $\left( \exists x\ \delta (x,\mathbf{u})\rightarrow \exists x<b\ \delta
(x,\mathbf{v})\right) $ holds in $\mathcal{M}$ for all $\Delta _{0}$%
-formulae $\delta (x,\mathbf{y}).$

\begin{itemize}
\item Lemma 6.1.5 below enables us to carry out a routine back-and-forth
construction to build sequences $\left( u_{k}:k<\omega \right) $ and $\left(
v_{k}:k<\omega \right) $ that satisfy (I), (II), and (III), thereby
establishing $(2)\Rightarrow (3)$ of Theorem 6.1. However, the proof of
Lemma 6.1.5 is labyrinthine, so we beg for the reader's indulgence.\medskip
\end{itemize}

\noindent \textbf{6.1.5.}$~$\textbf{Lemma.}$~$\textit{Suppose }$\mathbf{u}%
=\left\langle u_{r}:r<n\right\rangle $ \textit{and} $\ \mathbf{v}%
=\left\langle v_{r}:r<n\right\rangle $ \textit{are in} $\mathcal{M}$ \textit{%
with} $\max (\mathbf{v)}<b,$ \textit{and} \textit{both }$P(\mathbf{u},%
\mathbf{v})$ \textit{and} $Q(\mathbf{u},\mathbf{v})$ \textit{hold. Then}%
:\smallskip

\noindent \textbf{(a)} \textit{For every} $u^{\prime }\in M$ \textit{there
is }$v^{\prime }<b$ \textit{such that} \textit{both} $P(\left\langle \mathbf{%
u},u^{\prime }\right\rangle ,\left\langle \mathbf{v},v^{\prime
}\right\rangle )$ \textit{and} $Q(\left\langle \mathbf{u},u^{\prime
}\right\rangle ,\left\langle \mathbf{v},v^{\prime }\right\rangle )$ \textit{%
hold}; \textit{and}\smallskip

\noindent \textbf{(b)} \textit{For every} $v^{\prime }\in M$ \textit{with} $%
v^{\prime }<\max (\mathbf{v})$ \textit{there is some} $u^{\prime }\in M$
\textit{such that} \textit{both} $P(\left\langle \mathbf{u},u^{\prime
}\right\rangle ,\left\langle \mathbf{v},v^{\prime }\right\rangle )$ \textit{%
and} $Q(\left\langle \mathbf{u},u^{\prime }\right\rangle ,\left\langle
\mathbf{v},v^{\prime }\right\rangle )$ \textit{hold}.\medskip

\noindent \textbf{Proof of (a) of Lemma 6.1.5.}$~$We begin by noting that it
is well-known \cite[Lem.~2]{Yokoyama} that if $P(\mathbf{u},\mathbf{v})$
holds, then the proof of the basic Friedman embedding theorem as in \cite[%
Thm.~12.3]{Kaye's text} works for countable nonstandard models of $\mathrm{I}%
\Sigma _{1}$ and therefore:\smallskip

\noindent (1) There is a proper initial self-embedding $j_{0}$ of $\mathcal{M%
}$ such that $j_{0}(\mathcal{M})<b$ and $j_{0}(\mathbf{u})=\mathbf{v}$%
.\smallskip

\noindent Given $u^{\prime }\in M$ consider the type $p(x)=p_{1}(x)\cup
p_{2}(x)$, where:

\begin{center}
$p_{1}(x):=\{x<b\}\cup \{\left[ f(\mathbf{v},x)\downarrow \right] ^{<b}:$ $%
f\in \mathcal{F}\ \mathrm{and}\ \mathcal{M}\models \left[ f(\mathbf{u}%
,u^{\prime })\downarrow \right] \},$
\end{center}

\noindent and

\begin{center}
$p_{2}(x):=\left\{
\begin{array}{c}
\left[ f(\mathbf{v},x)\downarrow \right] ^{<b}\wedge f(\mathbf{v},x)\neq f(%
\mathbf{u},u^{\prime }):\smallskip \\
f\in \mathcal{F},\ \mathcal{M}\models \left[ f(\mathbf{u},u^{\prime
})\downarrow \right] \ \mathrm{and}\ f(\mathbf{u},u^{\prime })\notin K^{1}(%
\mathcal{M})%
\end{array}%
\right\} .$
\end{center}

\noindent Clearly if some $v^{\prime }$ realizes $p(x)$ in $\mathcal{M}$,
then $v^{\prime }<b$ and both $P(\left\langle \mathbf{u},u^{\prime
}\right\rangle ,\left\langle \mathbf{v},v^{\prime }\right\rangle )$ and $%
Q(\left\langle \mathbf{u},u^{\prime }\right\rangle ,\left\langle \mathbf{v}%
,v^{\prime }\right\rangle )$ hold. The fact that $\mathrm{Sat}_{\Sigma _{1}}$
is $\Sigma _{1}$-definable in $\mathcal{M}$, coupled with part (b) of
Theorem 2.3, makes it clear that $p_{1}(x)\in \mathrm{SSy}(\mathcal{M})$. To
show that $p_{2}(x)\in \mathrm{SSy}(\mathcal{M})$, let $\left\langle \delta
_{i}:i\in M\right\rangle $ be a canonical enumeration of $\Delta _{0}$%
-formulae within $\mathcal{M}$, and let $f_{i}$ be the partial recursive
function defined in $\mathcal{M}$ via:

\begin{center}
$f_{i}(x)=y$ iff $\mathrm{Sat}_{\Sigma _{1}}\left[ \exists z\ \left( z=\mu
t\ \delta _{i}\left( x,\left( t\right) _{0},\left( t\right) _{1}\right)
\right) \wedge y=\left( z\right) _{0}\right] .$
\end{center}

\noindent Consider the subset $R$ of $\mathbb{N}$ defined as follows:

\begin{center}
$R:=\{\left\langle i,j\right\rangle \in \mathbb{N}:\mathcal{M}\models \left[
f_{j}(\mathbf{u},u^{\prime })\downarrow \right] ,\ \left[ f_{i}(0)\downarrow %
\right] ,\ \mathrm{and}\ f_{j}(\mathbf{u},u^{\prime })=f_{i}(0)\}.$
\end{center}

\noindent Using the fact that $\mathrm{Sat}_{\Sigma _{1}}$ has a $\Sigma
_{1} $-description one can readily verify that $R$ is the intersection with $%
\mathbb{N}$ of a subset of $M$ that is parametrically $\Sigma _{1}$%
-definable in $\mathcal{M}$, so $R\in \mathrm{SSy}(\mathcal{M})$. Moreover,
using Lemma 6.1.4 we have:

\begin{center}
$\overset{A}{\overbrace{\{j\in \mathbb{N}:\exists i\ \left\langle
i,j\right\rangle \in R\}}}\ \medskip =$

$\overset{B}{\overbrace{\ \{j\in \mathbb{N}:\mathcal{M}\models \left[ f_{j}(%
\mathbf{u},u^{\prime })\downarrow \right] \ \mathrm{and}\ f_{j}(\mathbf{u}%
,u^{\prime })\in K^{1}(\mathcal{M})\}}.}$
\end{center}

\noindent Clearly $A$ is arithmetical in $R$, so $A\in \mathrm{SSy}(\mathcal{%
M})$ since we are assuming that $\mathbb{N}$ is strong in $\mathcal{M}$
(recall that by Theorem 2.5, $\mathrm{SSy}(\mathcal{M})$\ is arithmetically
closed). Hence $B\in \mathrm{SSy}(\mathcal{M}).$ Coupled with the closure of
$\mathrm{SSy}(\mathcal{M})$ under Turing reducibility and Boolean
operations, this shows that $p_{2}(x)\in \mathrm{SSy}(\mathcal{M}),$ which
finally makes it clear that $p(x)\in \mathrm{SSy}(\mathcal{M})$. \medskip

On the other hand, each formula in $p_{1}(x)$ is a $\Delta _{0}$-formula
(with parameters $v$ and $b)$, and each formula in $p_{2}(x)$ is a $\Sigma
_{1}$-formula (with parameters $\mathbf{u}$, $\mathbf{v}$, and $v^{\prime }$%
). In light of Remark 2.3.1, to show that $p(x)$ is realizable in $\mathcal{M%
}$ it is sufficient to verify that $p(x)$ is finitely realizable in $%
\mathcal{M}$. \medskip

Suppose $p(x)$ is not finitely realizable in $\mathcal{M}$. Note that the
formulae in $p_{1}(x)$ are closed under conjunctions, and that by (1) $%
p_{1}(x)$ is finitely realizable in $\mathcal{M}$. So for some $f\in
\mathcal{F}$, and some nonempty finite $\{g_{i}:i\leq k\}\subseteq $ $%
\mathcal{F}$ we have:\smallskip

\noindent (2) $\mathcal{M}\models \left[ f(\mathbf{u},u^{\prime })\downarrow %
\right] .$\smallskip

\noindent (3) $\mathcal{M}\models \left[ g_{i}(\mathbf{u},u^{\prime
})\downarrow \right] $ and $g_{i}\left( \mathbf{u},u^{\prime }\right) \notin
K^{1}(\mathcal{M})$ for $i\leq k.\footnote{%
As a warm-up, the reader may first wish to focus on the special but
instructive case $k=0$ in the argument that follows.}$\smallskip

\noindent (4) $\mathcal{M}\models \forall x<b\left(
\begin{array}{c}
\left[ f(\mathbf{v},x)\downarrow \right] ^{<b}\rightarrow \\
\bigvee\limits_{i=0}^{k}\left( \left[ g_{i}(\mathbf{v},x)\downarrow \right]
^{<b}\rightarrow g_{i}(\mathbf{v},x)=g_{i}\left( \mathbf{u},u^{\prime
}\right) \right)%
\end{array}%
\right) .\smallskip $

\noindent We may assume that $k$ is minimal in the sense that for any $%
f^{\prime }\in \mathcal{F}$ such that (2) holds with $f$ replaced by $%
f^{\prime }$ and any $k^{\prime }<k$, there is no subset $\{g_{i}^{\prime
}:i\leq k^{\prime }\}$ of $\mathcal{F}$ which has the property that both (3)
and (4) hold when $k$ is replaced by $k^{\prime }$, $f$ is replaced by $%
f^{\prime }$, and $g_{i}$ is replaced by $g_{i}^{\prime }.$ \medskip

\noindent By existentially quantifying $g_{0}(\mathbf{u},u^{\prime }),\cdot
\cdot \cdot ,g_{k}(\mathbf{u},u^{\prime })$ in (4) we obtain:\smallskip

\noindent (5) $\mathcal{M}\models \exists y$ $\theta (b,\mathbf{v,}y),$
where:

\begin{center}
$\theta (b,\mathbf{v},y):=\forall x<b\left(
\begin{array}{c}
\left[ f(\mathbf{v},x)\downarrow \right] ^{<b}\rightarrow \\
\bigvee\limits_{i=0}^{k}\left( \left[ g_{i}(\mathbf{v},x)\downarrow \right]
^{<b}\rightarrow g_{i}(\mathbf{v},x)=\left( y\right) _{i}\right)%
\end{array}%
\right) .$\smallskip
\end{center}

\noindent At this point we wish to define functions $h_{i}\in \mathcal{F}$
for $i\leq k$. We will denote the input of each $h_{i}$ by the symbol $%
\Diamond $ for better readability. For $i\leq k$, first let:

\begin{center}
$w_{0}(\Diamond ):=\mu w\ \exists y<w\ \theta (w,\Diamond ,y)$, and

$h(\Diamond ):=\mu y<w_{0}(\Diamond )$ $\theta (w_{0}(\Diamond ),\Diamond
,y),$
\end{center}

\noindent and then define:

\begin{center}
$h_{i}(\Diamond ):=\left( h(\Diamond )\right) _{i}.$
\end{center}

\noindent Clearly for each $i\leq k$, $h_{i}\in \mathcal{F}$; and $w_{0}$ is
well-defined iff $\left[ h_{i}(\Diamond )\downarrow \right] $ for each $%
i\leq k.$ The definition of $h_{i}$ together with (5) and the assumption
that $\max (\mathbf{v})<b$ makes it clear that:\smallskip

\noindent (6) $\mathcal{M}\models \varphi (b,\mathbf{v})$, where $\varphi (b,%
\mathbf{v})$ is the formula expressing\footnote{%
Note that a stronger form of statement (6) in which $\left[ h_{i}(\mathbf{v}%
)\downarrow \right] ^{<\max (\mathbf{v})+1}$ is weakened to $\left[ h_{i}(%
\mathbf{v})\downarrow \right] ^{<b}$ also holds, but (6) turns out to be the
appropriate ingredient for the argument that follows.}:

\begin{center}
$\forall x<b\left(
\begin{array}{c}
\left[ f(\mathbf{v},x)\downarrow \right] ^{<b}\rightarrow \\
\bigvee\limits_{i=0}^{k}\left( \left( \left[ g_{i}(\mathbf{v},x)\downarrow %
\right] ^{<b}\wedge \left[ h_{i}(\mathbf{v})\downarrow \right] ^{<\max (%
\mathbf{v})+1}\right) \rightarrow g_{i}(\mathbf{v},x)=h_{i}(\mathbf{v}%
)\right)%
\end{array}%
\right) .$
\end{center}

\noindent A salient feature of $\varphi (b,\mathbf{v})$ is that it is
expressible as a $\Pi _{1}^{<b}$-formula, i.e., a formula of the form $%
\forall z<b\ \delta (\mathbf{v},z),$ where $\delta $ is $\Delta _{0}.$
Recall that by assumption $P(\mathbf{u},\mathbf{v})$ holds, and that by
contraposition $P(\mathbf{u},\mathbf{v})$ is equivalent to:

\begin{center}
\textquotedblleft For all $\Delta _{0}$-formulae $\delta $, $\mathcal{M}%
\models \forall z<b\ \delta (\mathbf{v},z)\rightarrow \forall z\ \delta (%
\mathbf{u},z)$\textquotedblright .
\end{center}

\noindent So by (6) and $P(\mathbf{u},\mathbf{v})$, we may deduce: \smallskip

\noindent (7) $\mathcal{M}\models \forall x\left(
\begin{array}{c}
\left[ f(\mathbf{u},x)\downarrow \right] \rightarrow \\
\bigvee\limits_{i=0}^{k}\left( \left( \left[ g_{i}(\mathbf{u},x)\downarrow %
\right] \wedge \left[ h_{i}(\mathbf{u})\downarrow \right] ^{<\max (\mathbf{u}%
)+1}\right) \rightarrow g_{i}(\mathbf{u},x)=h_{i}(\mathbf{u})\right)%
\end{array}%
\right) ,$\smallskip

\noindent Recall that by (2) $\mathcal{M}\models \left[ f(\mathbf{u}%
,u^{\prime })\downarrow \right] $, so in light of (7) we have:\smallskip

\noindent (8) $\mathcal{M}\models \bigvee\limits_{i=0}^{k}\left( \left( %
\left[ g_{i}(\mathbf{u},u^{\prime })\downarrow \right] \wedge \left[ h_{i}(%
\mathbf{u})\downarrow \right] ^{<\max (\mathbf{u})+1}\right) \rightarrow
g_{i}(\mathbf{u},u^{\prime })=h_{i}(\mathbf{u})\right) .$\smallskip

\noindent Based on (8) we may assume without loss of generality:\smallskip

\noindent (9) $\mathcal{M}\models \left( \left[ g_{0}(\mathbf{u},u^{\prime
})\downarrow \right] \wedge \left[ h_{0}(\mathbf{u})\downarrow \right]
^{<\max (\mathbf{u})+1}\right) \rightarrow g_{0}(\mathbf{u},u^{\prime
})=h_{0}(\mathbf{u}).$\smallskip

\noindent At this point we claim that the following statement $(\ast )$\ is
true. Note that since the subformula marked as $\psi $ in $(\ast )$ (the
premise of the implication) is equivalent to a formula in $p_{1}(x)$ and the
index $i$ in the disjunction in $(\ast )$ starts from $i=1$, the veracity of
$(\ast )$ contradicts the minimality of $k$.\smallskip

\noindent $(\ast )$ $\mathcal{M}\models $\smallskip $\forall x<b\left(
\begin{array}{c}
\overset{\psi }{\overbrace{\left(
\begin{array}{c}
\left[ f(\mathbf{v},x)\downarrow \right] ^{<b}\wedge \left[ g_{0}(\mathbf{v}%
,x)\downarrow \right] ^{<b}\wedge \\
\left[ h_{0}(\mathbf{v})\downarrow \right] ^{<\max (\mathbf{v})+1}\wedge
g_{0}(\mathbf{v},x)=h_{0}(\mathbf{v)}%
\end{array}%
\right) }}\rightarrow \smallskip \\
\left( \bigvee\limits_{i=1}^{k}\left[ g_{i}(\mathbf{v},x)\downarrow \right]
^{<b}\rightarrow g_{i}\left( \mathbf{v},x\right) =g_{i}\left( \mathbf{u}%
,u^{\prime }\right) \right)%
\end{array}%
\right) .$\smallskip

\noindent Suppose to the contrary that $(\ast )$ fails. Then for some $c\in
M $:\smallskip

\noindent (10) $\mathcal{M}\models $\smallskip $\left( c<b\right) \wedge
\left(
\begin{array}{c}
\left(
\begin{array}{c}
\left[ f(\mathbf{v},c)\downarrow \right] ^{<b}\wedge \left[ g_{0}(\mathbf{v}%
,c)\downarrow \right] ^{<b}\wedge \\
\left[ h_{0}(\mathbf{v})\downarrow \right] ^{<\max (\mathbf{v})+1}\wedge
g_{0}(\mathbf{v},c)=h_{0}(\mathbf{v)}%
\end{array}%
\right) \wedge \smallskip \\
\lnot \left( \bigvee\limits_{i=1}^{k}\left[ g_{i}(\mathbf{v},c)\downarrow %
\right] ^{<b}\rightarrow g_{i}\left( \mathbf{v},c\right) =g_{i}\left(
\mathbf{u},u^{\prime }\right) \right)%
\end{array}%
\right) .$\smallskip

\noindent Recall that by (3) $\left[ g_{0}(\mathbf{u},u^{\prime })\downarrow %
\right] $, which coupled with (10) makes it clear that: $\smallskip $

\noindent (11) $\mathcal{M}\models \left[ g_{0}(\mathbf{u},u^{\prime
})\downarrow \right] \wedge \left[ g_{0}(\mathbf{v},c)\downarrow \right]
\wedge \left[ h_{0}(\mathbf{v})\downarrow \right] ^{<\max (\mathbf{v})+1}.$%
\smallskip

\noindent In light of (11), (10), and (4) we also have:\smallskip

\noindent (12) $\mathcal{M}\models h_{0}(\mathbf{v})=g_{0}\left( \mathbf{v}%
,c\right) =g_{0}\left( \mathbf{u},u^{\prime }\right) .$\smallskip

\noindent By (12) $\mathcal{M}\models g_{0}\left( \mathbf{u},u^{\prime
}\right) =h_{0}(\mathbf{v)},$ and by (3) $g_{0}\left( \mathbf{u},u^{\prime
}\right) \notin K^{1}(\mathcal{M})$, hence:\smallskip

\noindent (13) $h_{0}(\mathbf{v})\notin K^{1}(\mathcal{M}).$\smallskip

\noindent On the other hand, by (11) $\mathcal{M}\models \left[ h_{0}(%
\mathbf{v})\downarrow \right] ^{<\max (\mathbf{v})+1}$, so (1) makes it
clear that $\mathcal{M}\models \left[ h_{0}(\mathbf{u})\downarrow \right]
^{<\max (\mathbf{u})+1}$. Therefore in light of (11) and (9) $g_{0}(\mathbf{u%
},u^{\prime })=h_{0}(\mathbf{u}),$ so $h_{0}(\mathbf{u})=h_{0}(\mathbf{v})$,
which by our assumption that $Q(\mathbf{u},\mathbf{v})$ holds, implies $%
h_{0}(\mathbf{v})\in K^{1}(\mathcal{M})$, thereby contradicting (13). This
contradiction demonstrates that $(\ast )$ is true, thus refuting the
minimality of $k$ and completing the proof.\hfill $\square $ Lemma
6.1.5(a)\medskip

\noindent \textbf{Proof of (b) of Lemma 6.1.5.}$~$The proof of this part has
some resemblances to the proof of part (a), but it also exhibits certain
differences. Let $\max (\mathbf{v})=v_{j}.$ Then by the assumption that $P(%
\mathbf{u},\mathbf{v})$ holds, $\max (\mathbf{u})=u_{j}$. Given $v^{\prime
}\in M$ with $v^{\prime }<v_{j}$ consider the following type $%
q(x)=q_{1}(x)\cup q_{2}(x)$, where:

\begin{center}
$q_{1}(x):=\{x<u_{j}\}\cup \{\lnot \left[ f(\mathbf{u},x)\downarrow \right]
: $ $f\in \mathcal{F},$ $\mathcal{M}\models \lnot \left[ f(\mathbf{v}%
,v^{\prime })\downarrow \right] ^{<b}\},$
\end{center}

\noindent and

\begin{center}
$q_{2}(x):=\left\{
\begin{array}{c}
\left[ f(\mathbf{u},x)\downarrow \right] \rightarrow f(\mathbf{v},v^{\prime
})\neq f(\mathbf{u},x):\smallskip \\
f\in \mathcal{F},\ \left[ f(\mathbf{v},v^{\prime })\downarrow \right]
^{<b},\ \mathrm{and}\ f\left( \mathbf{v},v^{\prime }\right) \notin K^{1}(%
\mathcal{M})%
\end{array}%
\right\} .$
\end{center}

\noindent It is routine to verify that if some element $u^{\prime }$ of $%
\mathcal{M}$ realizes $q(x)$, then both $P(\left\langle \mathbf{u},u^{\prime
}\right\rangle ,\left\langle \mathbf{v},v^{\prime }\right\rangle )$ and $%
Q(\left\langle \mathbf{u},u^{\prime }\right\rangle ,\left\langle \mathbf{v}%
,v^{\prime }\right\rangle )$ hold. Also one can show that $q(x)\in \mathrm{%
SSy}(\mathcal{M})$ using a reasoning analogous to the one used in the proof
of part (a) to show that $p(x)\in \mathrm{SSy}(\mathcal{M}).$ By Remark
2.3.1 to show that $q(x)$ is realized in $\mathcal{M}$ it suffices to
demonstrate that $q(x)$ is finitely realizable in $\mathcal{M}$ since $q(x)$
is a short $\Pi _{1}$-type . Suppose $q(x)$ is not finitely realized in $%
\mathcal{M}$. Then since the formulae in $q_{1}(x)$ are closed under
conjunctions and $q_{1}(x)$ is finitely satisfiable in $\mathcal{M}$ by
statement (1) of the proof of Lemma 6.1.5(a), for some $f\in \mathcal{F}$,
and some nonempty finite $\{g_{i}:i\leq k\}\subseteq \mathcal{F}$, where $k$
is minimal, we have:

\noindent (1) $\mathcal{M}\models \lnot \left[ f(\mathbf{v},v^{\prime
})\downarrow \right] ^{<b}.$\smallskip

\noindent (2) $\mathcal{M}\models \left[ g_{i}(\mathbf{v},v^{\prime
})\downarrow \right] ^{<b}$ and $g_{i}\left( \mathbf{v},v^{\prime }\right)
\notin K^{1}(\mathcal{M})$ for $i\leq k.$\smallskip

\noindent (3) $\mathcal{M}\models \forall x<u_{j}\left(
\begin{array}{c}
\left[ f\left( \mathbf{u},x\right) \downarrow \right] \vee \\
\bigvee\limits_{i=0}^{k}\left( \left[ g_{i}(\mathbf{u},x)\downarrow \right]
\wedge g_{i}(\mathbf{u},x)=g_{i}\left( \mathbf{v},v^{\prime }\right) \right)%
\end{array}%
\right) .$\smallskip

\noindent By existentially quantifying $g_{0}\left( \mathbf{v},v^{\prime
}\right) ,\cdot \cdot \cdot ,g_{k}\left( \mathbf{v},v^{\prime }\right) $ in
(4), and taking advantage of the veracity of $\mathrm{B}\Sigma _{1}$ in $%
\mathcal{M}$ we obtain:\smallskip

\noindent (4) $\mathcal{M}\models \exists w\exists y<w\ \theta (w,\mathbf{u}%
,y),$ where:

\begin{center}
$\theta (w,\mathbf{u},y):=\forall x<u_{j}\left(
\begin{array}{c}
\left[ f\left( \mathbf{u},x\right) \downarrow \right] ^{<w}\vee \\
\bigvee\limits_{i=0}^{k}\left( \left[ g_{i}(\mathbf{u},x)\downarrow \right]
^{<w}\wedge g_{i}(\mathbf{u},x)=\left( y\right) _{i}\right)%
\end{array}%
\right) .$\smallskip
\end{center}

\noindent As in the proof of part(a), we will define functions $h_{i}\in
\mathcal{F}$ for $i\leq k$ and will denote the input of each $h_{i}$ by the
symbol $\Diamond .$ For $i\leq k,$ first define:

\begin{center}
$w_{0}(\Diamond ):=\mu w\ \exists y<w\ \theta (w,\mathbf{\Diamond },y),$ and

$h(\Diamond ):=\mu y<w_{0}(\Diamond )$ $\theta (w_{0}(\Diamond ),\Diamond
,y),$
\end{center}

\noindent and then define:

\begin{center}
$h_{i}(\Diamond ):=\left( h(\Diamond )\right) _{i}.$
\end{center}

\noindent Clearly $h_{i}\in \mathcal{F}$ for $i\leq k$; and $w_{0}$ is
well-defined iff $\left[ h_{i}(\Diamond )\downarrow \right] $ for each $%
i\leq k.$ The definition of $h_{i}$ together with (4) yields:\smallskip

\noindent (5) $\mathcal{M}\models \exists w\ \varphi (\mathbf{u})^{<w},$
where:

\begin{center}
$\varphi (\mathbf{u}):=\forall x<u_{j}\left(
\begin{array}{c}
\left[ f\left( \mathbf{u},x\right) \downarrow \right] \vee \\
\bigvee\limits_{i=0}^{k}\left( \left[ g_{i}(\mathbf{u},x)\downarrow \right]
\wedge \left[ h_{i}(\mathbf{u})\downarrow \right] \wedge g_{i}(\mathbf{u}%
,x)=h(\mathbf{u})\right)%
\end{array}%
\right) $
\end{center}

\noindent where $\varphi (\mathbf{u})^{<w}$ is the $\Delta _{0}$-formula
obtained by relativizing $\varphi (\mathbf{u})$ to the predecessors of $w$
(formally: the result of replacing every unbounded quantifier $\mathrm{Q}z$
in $\varphi (\mathbf{u})$ to $\mathrm{Q}z<w$). Also note that $\varphi (%
\mathbf{u})$ can be written as a $\Sigma _{1}$-formula since $\mathcal{M}%
\models \mathrm{B}\Sigma _{1}$. Therefore $\mathcal{M}\models \varphi (%
\mathbf{v})^{<b}$ by putting (5) together with our assumption that $P(%
\mathbf{u},\mathbf{v})$ holds, in other words we now have: \smallskip

\noindent (6) $\mathcal{M}\models \forall x<v_{j}\left(
\begin{array}{c}
\left[ f\left( \mathbf{v},x\right) \downarrow \right] ^{<b}\vee \\
\bigvee\limits_{i=0}^{k}\left( \left[ g_{i}(\mathbf{v},x)\downarrow \right]
^{<b}\wedge \left[ h_{i}(\mathbf{v})\downarrow \right] ^{<b}\wedge g_{i}(%
\mathbf{v},x)=h_{i}(\mathbf{v})\right)%
\end{array}%
\right) .$\smallskip

\noindent Putting (1) together with (6) and the assumption that $v^{\prime
}\leq \max (\mathbf{v})=v_{j}$ gives us: \smallskip

\noindent (7) $\mathcal{M}\models \bigvee\limits_{i=0}^{k}\left( \left[
g_{i}(\mathbf{v},v^{\prime })\downarrow \right] ^{<b}\wedge \left[ h_{i}(%
\mathbf{v})\downarrow \right] ^{<b}\wedge g_{i}(\mathbf{v},v^{\prime
})=h_{i}(\mathbf{v})\right) .$\smallskip

\noindent Based on (7) we may assume without loss of generality:\smallskip

\noindent (8) $\mathcal{M}\models \left( \left[ g_{0}(\mathbf{v},v^{\prime
})\downarrow \right] \wedge \left[ h_{0}(\mathbf{v})\downarrow \right]
^{<b}\wedge g_{0}(\mathbf{v},v^{\prime })=h_{0}(\mathbf{v})\right) .$%
\smallskip

\noindent At this point we claim that $(\ast )$\ below holds. Note that $%
(\ast )$ contradicts the minimality of $k$ since the subformula marked as $%
\psi $ in $(\ast )$ (embraced by curly braces) is equivalent to the negation
of a formula in $q_{1}(x)$, and in the disjunction in $(\ast )$ the index $i$
starts from $i=1$.\smallskip

\noindent $(\ast )$ $\mathcal{M}\models $\smallskip $\forall x<u_{j}\left(
\begin{array}{c}
\overset{\psi }{\overbrace{\left\{ \left[ f\left( \mathbf{u},x\right)
\downarrow \right] \vee \left(
\begin{array}{c}
\left( \left[ g_{0}(\mathbf{u},x)\downarrow \right] \wedge \left[ h_{0}(%
\mathbf{u})\downarrow \right] \right) \wedge \\
g_{0}(\mathbf{u},x)\neq h_{0}(\mathbf{u})%
\end{array}%
\right) \right\} }}\vee \\
\bigvee\limits_{i=1}^{k}\left( \left[ g_{i}(\mathbf{u},x)\downarrow \right]
\wedge \left[ h_{i}(\mathbf{u})\downarrow \right] \wedge g_{i}(\mathbf{u}%
,x)=g_{i}\left( \mathbf{v},v^{\prime }\right) \right)%
\end{array}%
\right) .$\smallskip

\noindent Suppose to the contrary that $(\ast )$ fails. Then for some $%
c<u_{j}$:\smallskip

\noindent (9) $\mathcal{M}\models $\smallskip $\left(
\begin{array}{c}
\left\{ \lnot \left[ f\left( \mathbf{u},c\right) \downarrow \right] \wedge
\left(
\begin{array}{c}
\left( \left[ g_{0}(\mathbf{u},c)\downarrow \right] \wedge \left[ h_{0}(%
\mathbf{u})\downarrow \right] \right) \rightarrow \\
g_{0}(\mathbf{u},c)=h_{0}(\mathbf{u})%
\end{array}%
\right) \right\} \wedge \\
\lnot \left( \bigvee\limits_{i=1}^{k}\left[ g_{i}(\mathbf{u},c)\downarrow %
\right] \wedge \left[ h_{i}(\mathbf{u})\downarrow \right] \wedge g_{i}(%
\mathbf{u},c)=g_{i}\left( \mathbf{v},v^{\prime }\right) \right)%
\end{array}%
\right) .$\smallskip

\noindent Recall that by (2) $\left[ g_{0}(\mathbf{v},v^{\prime })\downarrow %
\right] $. Since $c<u_{j}$, by putting (9) together with (3) we can conclude
that:\smallskip

\noindent (10) $\mathcal{M}\models \left[ g_{0}(\mathbf{u},c)\downarrow %
\right] \wedge g_{0}\left( \mathbf{u},c\right) =g_{0}\left( \mathbf{v}%
,v^{\prime }\right) .$\smallskip

\noindent Also, $\mathcal{M}\models \left[ h_{0}(\mathbf{u})\downarrow %
\right] $ by (5) and (9). So in light of (9) and (10) we have:\smallskip

\noindent (11) $\mathcal{M}\models g_{0}\left( \mathbf{v},v^{\prime }\right)
=g_{0}\left( \mathbf{u},c\right) =h_{0}(\mathbf{u}).$\smallskip

\noindent By (11) $g_{0}\left( \mathbf{v},v^{\prime }\right) =h_{0}(\mathbf{u%
}).$ So $h_{0}(\mathbf{u})\notin K^{1}(\mathcal{M})$ since $g_{0}\left(
\mathbf{v},v^{\prime }\right) \notin K^{1}(\mathcal{M})$ by (2). But then we
have a contradiction since by (8) $g_{0}\left( \mathbf{v},v^{\prime }\right)
=h_{0}(\mathbf{v}),$ hence $h_{0}(\mathbf{v})=h_{0}(\mathbf{u})$, so $h_{0}(%
\mathbf{u})\in K^{1}(\mathcal{M})$ by the assumption that $Q(\mathbf{u},%
\mathbf{v})$ holds. This concludes our proof of $(\ast ),$ which in turn
contradicts the minimality of $k$ and finishes the proof. \hfill $\square $
Lemma 6.1.5(b)\medskip

\noindent With Lemma 6.1.5 at hand, the proof of $(2)\Rightarrow (3)$ of
Theorem 6.1 is now complete.\hfill $\square $

\textbf{\bigskip }

\begin{center}
\textbf{7.~CLOSING REMARKS AND OPEN QUESTIONS\bigskip }
\end{center}

\noindent \textbf{7.1.}$~$\textbf{Remark.}$~$Let $\mathcal{L}$ be a finite
extension of $\mathcal{L}_{A}$. An inspection of the proofs of Theorems 4.1,
5.1, and 6.1 make it clear that the equivalence of conditions (2) and (3) of
these theorems stays valid for countable nonstandard models of \textrm{I}$%
\Sigma _{1}(\mathcal{L})$. Furthermore, condition (1) of the aforementioned
theorems remains equivalent to the other two conditions in the setting of
\textrm{I}$\Sigma _{1}(\mathcal{L})$ if (1) is strengthened to the assertion
that $j$ is a $\Delta _{0}(\mathcal{L})$-elementary self-embedding of $%
\mathcal{M}$.\medskip

\noindent \textbf{7.2.}$~$\textbf{Remark.}$~$Wilkie \cite{Wilke thesis}
showed that if $\mathcal{M}$ is a countable nonstandard model of $\mathrm{PA}
$, then:

\begin{center}
$\left\vert \{I:I\ \mathrm{is\ a\ cut\ of\ }\mathcal{M}\ \mathrm{and}\
I\cong \mathcal{M}\}\right\vert =2^{\aleph _{0}}$.
\end{center}

\noindent The proof strategy in \cite[Thm.~2.7 ($n=0$)]{Kaye's text} of
Wilkie's theorem can be shown to work for all countable nonstandard models $%
\mathcal{M}$ of $\mathrm{I}\Sigma _{1}.$ Moreover, Theorem 4.1 can be
refined by strengthening condition (3) of that theorem to state that there
are\textit{\ }$2^{\aleph _{0}}$\textit{-}many cuts of $\mathcal{M}$ that can
appear as the range of initial embeddings $j$ of $\mathcal{M}$ for which $I=%
\mathrm{I}_{\mathrm{fix}(j)}.$ These results will appear in \cite{Bahrami
Diss.}. \medskip

\noindent \textbf{7.3.}$~$\textbf{Remark.}$~$The main results of the paper
(Theorems 4.1, 5.1, and 6.1) lend themselves to a hierarchical
generalization in which $\mathcal{M}\models \mathrm{I}\Sigma _{n+1}$ and the
self-embedding $j$ is stipulated to be $\Sigma _{n}$-elementary\textit{. }%
These results will also appear in \cite{Bahrami Diss.}.\medskip

\noindent \textbf{7.4.}$~$\textbf{Question.}$~$\textit{Is it true that in
Theorems 5.1 and 6.1 condition (3) can be strengthened by adding that there
are continuum-many cuts of }$\mathcal{M}$\textit{\ that can be realized as
the range of }$j$\textit{?}

\noindent Remark 7.1 suggests that Question 7.4 has a positive
answer.\medskip

\noindent \textbf{7.5.}$~$\textbf{Question.}$~$\textit{Is there some} $n\in
\omega $ \textit{such that} \textit{every countable nonstandard model of }%
\textrm{I}$\Sigma _{n}$\textit{\ has a contractive }(\textit{i.e., }$%
j(a)\leq a$\ \textit{always})\ \textit{proper initial self-embedding}?
\textit{And if the answer is positive, what is the minimal such} $n$?

\noindent The above question is motivated by Corollary 3.4.2.\medskip

\noindent \textbf{7.6.}$~$\textbf{Question.}$~$\textit{Suppose }$\mathcal{M}$%
\textit{\ is a countable nonstandard model of }\textrm{I}$\Sigma _{1}$%
\textit{\ in which }$\mathbb{N}$\ \textit{is a strong cut}. \textit{Is every
proper }$\Sigma _{1}$\textit{-elementary submodel of $\mathcal{M}$
isomorphic to }$\mathrm{Fix}(j)$\textit{\ for some }$j$\textit{\ }$\in
\mathrm{PISE}(\mathcal{M)}$?

\noindent The above question is prompted by the result mentioned in footnote
3, and the fact that the proof of Theorem 6.1 makes it clear that the
theorem remains valid if in conditions and (1) and (3) of the statements of
that theorem, the requirement that $\mathrm{Fix}(j)=K^{1}(\mathcal{M})$ is
modified to $\mathrm{Fix}(j)=K^{1}(\mathcal{M},m)$, where $m\in M.$\medskip

\noindent \textbf{7.7.}$~$\textbf{Question.}$~$\textit{Suppose} $I$ \textit{%
is a strong cut of} $\mathcal{M}\models \mathrm{I}\Sigma _{1}$, $\mathcal{N}%
\prec _{\Sigma _{1}}\mathcal{M}$\textit{, and }$\mathcal{N}$ \textit{is} $I$%
\textit{-coded} (\textit{i.e., there is an element} $s$ \textit{of} $%
\mathcal{M}$ \textit{such that} $N=\left\{ (s)_{i}:i\in I\right\} $ \textit{%
and} $s_{i}\neq s_{j}$ \textit{if} $i<j\in I)$\textit{, then} $\mathcal{N}$\
\textit{can be realized as} $\mathrm{Fix}(j)$ \textit{for some} $j\in
\mathrm{PISE}(\mathcal{M)}$?

\noindent The impetus for the above question can be found in \cite[Thm.~4.5.1%
]{Me unified}. 

\begin{center}
\begin{tabular}{ll}
{\small Saeideh Bahrami} & {\small Ali Enayat} \\
{\small Dept.~of Mathematics} & {\small Dept.~of Philosophy, Linguistics, \&
Theory of Science} \\
{\small Tarbiat Modares University} & {\small University of Gothenburg} \\
{\small P.O.~Box 14115-111, Tehran, Iran} & {\small Box 200, SE 405 30,
Gothenburg, Sweden} \\
{\small E-mail: bahrami.saeideh@gmail.com} & {\small E-mail: ali.enayat@gu.se%
}%
\end{tabular}
\end{center}

\end{document}